\documentclass{article}
\usepackage[T2A]{fontenc}
\usepackage[cp1251]{inputenc}
\usepackage[english,russian]{babel}
\usepackage[tbtags]{amsmath}
\usepackage{amsfonts,amssymb,mathrsfs,amscd}
\numberwithin{equation}{section}
\newtheorem{remark}{Remark}[section]
\newtheorem{lemma}{Lemma}[section]
\newtheorem{theorem}{Theorem}[section]
\newtheorem{definition}{Definition}[section]
\newtheorem{corollary}{Corollary}[section]

\renewcommand{\Im}{\mathop{\rm Im}}
\renewcommand{\Re}{\mathop{\rm Re}}

\newcommand{\supp}{\mathop{\rm supp}}

\begin{document}
\begin{Large}
\thispagestyle{empty}
\begin{center}
{\bf Inverse scattering problem for a third-order operator with non-local potential\\
\vspace{5mm}
V. A. Zolotarev}\\

B. Verkin Institute for Low Temperature Physics and Engineering
of the National Academy of Sciences of Ukraine\\
47 Nauky Ave., Kharkiv, 61103, Ukraine

Department of Higher Mathematics and Informatics, V. N. Karazin Kharkov National University \\
4 Svobody Sq, Kharkov, 61077,  Ukraine

\end{center}
\vspace{5mm}

{\small {\bf Abstract.} Direct and inverse scattering problem for an operator with non-local potential is solved in the paper. The method is based on the Riemann boundary value problem on a bundle of three straight lines. Description of scattering problem data is given.}
\vspace{5mm}

{\it Mathematics
Subject Classification 2020:} 34L10, 34L15.\\

{\it Key words}: inverse spectral problem, boundary value problem, Jost solutions.
\vspace{5mm}

\begin{center}
{\bf Introduction}
\end{center}
\vspace{5mm}

Inverse scattering problem is one of intriguing problems of mathematical physics. I. M. Gelfand -- B. M. Levitan, V. A. Marchenko, M. G. Krein in their works \cite{1} -- \cite{5} were the first to solve it for the one-dimensional Shr$\ddot{\rm o}$dinger equation. Method of inverse scattering problem became an effective tool for integrating nonlinear partial differential equations (Korteweg -- De Vries equation, sine -- Gordon equation etc.) \cite{6}. First, a non-linear equation is matched with L -- A Lax pair \cite{6} where L, as a rule, is a well-known second-order operator (e. g., Sturm -- Liouville operator) and then the method of inverse problem is applied to the operator L.
Search for L -- A pairs for Camassa -- Holm and Degasperis -- Procesi equations (waves in dispersive media) leads to L -- A pairs in which the operator L is of third order (cubic string) \cite{7} -- \cite{12}. This is the reason for writing this paper. It studies scattering scheme for a third-order differential operator on half-axis with non-local potential. Study of differential operators  of an order higher than two has its own specific features since in this case transformation operators that play the  key role for Sturm -- Liouville operators \cite{1, 2} are absent. Scattering problem for a second order differential operator with non-local potential is studied in \cite{13}, see also \cite{14, 15}.

Integro-differential operator
\begin{equation}
(L_\alpha y)(x)=iy'''(x)+\alpha\int\limits_0^\infty y(t)\overline{q}(t)dtq(x)\quad(x\in\mathbb{R}_+)\label{eq0.1}
\end{equation}
is the main object under study in this paper, here $\alpha\in\mathbb{R}$, $q\in L^2(\mathbb{R}_+)$. This operator is a one-dimensional perturbation of the operator $iD^3=P^3$ (${\displaystyle D=\frac d{dx}}$, ${\displaystyle P=\frac 1iD}$) coinciding with the third power of momentum operator $P$. Operator $L_\alpha$ \eqref{eq0.1} is symmetric under the boundary conditions $y(0)=0$, $y'(0)=0$ (for example) and does not have a self-adjoint expansion in this space \cite{16, 17, 18}.  And since spectral problem is based upon spectral analysis of a self-adjoint operator \cite{19, 3, 16, 17}, it is natural to complete the operator $L_\alpha$ to self-adjoint with departure from the original space. Unlike standard expansions \cite{16, 17, 18}, in this work a simpler (but not unique) expansion of the operator $L_\alpha$ is given.

Fourier transform upon $\mathbb{R}$ is closely related to the operator of second derivative and is constructed by two exponents (fundamental system of two oscillating solutions to a second order equation) to the one of which there corresponds the right half-axis, and to the other, correspondingly, the left. Section 1 outlines the basics of Fourier transform corresponding to a third-order equation (here fundamental system of equations is formed by three exponents: oscillating, decreasing and increasing). Therefore Fourier transform in this case is constructed on the system (bundle) of three straight lines in the plane $\mathbb{C}$ at the angle $2\pi/3$ between each other. In this section, analogues of the cosine and sine Fourier transforms are obtained.

Section 2 constructs Jost solutions (three of them) for a self-adjoint expansion of the operator $L_{\alpha}$ \eqref{eq0.1} and describes the set of such values of $\lambda$ for which these Jost solutions are linearly independent.

The third section of the work studies scattering problem that has its own specific features. Here, there is a single incident wave and two scattered ones and thus two scattering coefficients. This section describes these scattering coefficients and shows that they set a Riemann boundary value problem on a complex contour formed by three lines (Sec. 1). Similar boundary value problems arose in the work \cite{20}.

Section 4 gives method for solving inverse scattering problem based on a Riemann boundary value problem. It is shown that, under normalization condition, number $\alpha$ and function $q(x)$ are defined by scattering data.  Description of scattering data for this inverse problem is given.

\section{Preliminary information. Fourier transform on a system of straight lines}\label{s1}

{\bf 1.1.} Equation
\begin{equation}
iD^3y(x)=\lambda^3y(x)\quad(x\in\mathbb{R},D=\frac d{dx},\lambda\in\mathbb{C})\label{eq1.1}
\end{equation}
has \cite{21, 22} three linearly independent solutions $\{e^{i\lambda\zeta_kx}\}_1^3$, here $\{\zeta_k\}_1^3$ are roots of a cubic equation $\zeta^3=1$,
\begin{equation}
\zeta_1=1,\quad\zeta_2=\frac12(-1+i\sqrt3),\quad\zeta_3=\frac12(-1-i\sqrt3),\label{eq1.2}
\end{equation}
besides,
\begin{equation}
\zeta_1+\zeta_2+\zeta_3=0,\,\zeta_1\zeta_2\zeta_3=1,\,\zeta_1-\zeta_2=i\sqrt3\zeta_3,\,\zeta_2-\zeta_3=i\sqrt3\zeta_1,\,\zeta_3-\zeta_1=i\sqrt3\zeta_2.\label{eq1.3}
\end{equation}
Any solution to equation \eqref{eq1.1} is a linear combination of exponents $\{e^{i\lambda\zeta_kx}\}_1^3$. Instead of $\{e^{i\lambda\zeta_kx}\}_1^3$ one can take another system of fundamental solutions to equation \eqref{eq1.1}, for example, $s_0(i\lambda x)$, $s_1(i\lambda x)$, $s_2(i\lambda x)$, here
\begin{equation}
s_0(z)=\frac13\sum\limits_ke^{z\zeta_k};\quad s_1(z)=\frac13\sum\limits_k\frac1{\zeta_k}e^{z\zeta_k};\quad s_2(z)=\frac13\sum\limits_k\frac1{\zeta_k^2}e^{z\zeta_k}\label{eq1.4}
\end{equation}
($z\in\mathbb{C}$). The functions \eqref{eq1.4} are similar to cosines and sines for a second order equation ($-D^2y(x)=\lambda^2y(x)$). Their main properties are contained in the following statement \cite{21}.

\begin{lemma}\label{l1.1}
Entire functions of exponential type \eqref{eq1.4} meet the conditions

${\rm(i)}$ $s'_0(z)=s_2(z)$, $s'_1(z)=s_0(z)$, $s'_2(z)=s_1(z)$;

${\rm(ii)}$ $\overline{s_p(z)}=s_p(\overline{z})$, ($0\leq p\leq2$);

${\rm(iii)}$ $s_p(z\zeta_2)=\zeta_2^ps_p(z\zeta_2)$ ($0\leq p\leq2$);

${\rm(iv)}$ the Euler formula
$$e^{z\zeta_p}=s_0(z)+\zeta_ps_1(z)+\zeta_p^2s_2(z)\quad(1\leq p\leq3);$$

${\rm(v)}$ the functions \eqref{eq1.4} are solutions to the equation ${\displaystyle\frac{d^3}{dz^3}y(z)=y(z)}$ and satisfy the initial data
$$s_0(0)=1,\quad s'_0(0)=0,\quad s''_0(0)=1;$$
$$s_1(0)=0,\quad s'_1(0)=1,\quad s''_1(0)=0;$$
$$s_2(0)=0,\quad s'_2(0)=0,\quad s''_2(0)=1;$$

${\rm(vi)}$ the main identity
$$s_0^3(z)+s_1^3(z)+s_2^3(z)-3s_0(z)s_1(z)s_2(z)=1;$$

${\rm(vii)}$ addition formulas
$$s_0(z+w)=s_0(z)s_0(w)+s_1(z)s_2(w)+s_2(z)s_1(w);$$
$$s_1(z+v)=s_0(z)s_1(w)+s_1(z)s_0(w)+s_2(z)s_2(w);$$
\begin{equation}
s_2(z+w)=s_0(z)s_2(w)+s_1(z)s_1(w)+s_2(z)s_0(w);\label{eq1.5}
\end{equation}

${\rm{(viii)}}$ $3s_0(z)s_0(w)=s_0(z+w)+s_0(z+\zeta_2w)+s_0(z+\zeta_3w);$
$$3s_1(z)s_2(w)=s_0(z+w)+\zeta_2s_0(z+\zeta_2w)+\zeta_3s_0(z+\zeta_3w);$$
$$3s_0(z)s_1(w)=s_1(z+w)+s_1(z+\zeta_2w)+s_1(z+\zeta_3w);$$
$$3s_2(z)s_2(w)=s_1(z+w)+\zeta_2s_1(z+\zeta_2w)+\zeta_3s_1(z+\zeta_3w);$$
$$3s_2(z)s_0(w)=s_2(z+w)+s_2(z+\zeta_2w)+s_2(z+\zeta_3w);$$
$$3s_1(z)s_1(w)=s_2(z+w)+\zeta_3s_2(z+\zeta_2w)+\zeta_2s_2(z+\zeta_3w);$$

$\rm{(ix)}$ $3s_0^2(z)=s_0(2z)+2s_0(-z);$

$\hspace{10mm}3s_1^2(z)=s_2(2z)-2s_2(-z);$

$\hspace{10mm}3s_2^2(z)=s_1(2z)-2s_1(-z);$

$\rm{(x)}$ $s_0^2(z)-s_1(z)s_2(z)=s_0(-z);$

$\hspace{10mm}s_1^2(z)-s_2(z)s_0(z)=s_2(-z);$

$\hspace{10mm}s_2^2(z)-s_1(z)s_0(z)=s_1(-z);$

$\rm{(xi)}$ the Taylor formula

$\hspace{10mm}s_0(z)=1+{\displaystyle\frac{z^3}{3!}+\frac{z^6}{6!}+...}$

$\hspace{10mm}s_1(z)=z+{\displaystyle\frac{z^4}{4!}+\frac{z^7}{7!}+...}$

$\hspace{10mm}s_2(z)={\displaystyle\frac{z^2}2+\frac{z^5}{5!}+\frac{z^8}{8!}+...}.$
\end{lemma}

Proof of the lemma follows from \eqref{eq1.3}.

\begin{remark}\label{r1.1}
For equation ${\displaystyle\left(\frac1iD\right)^ny=\lambda^ny}$ for arbitrary $n\in\mathbb{N}$, similarly to \eqref{eq1.4}, it is easy to construct such fundamental system of solutions for which statements of Lemma \ref{l1.1} hold (of course, in the corresponding formulation).
\end{remark}

Solution to the Cauchy problem
\begin{equation}
iD^3y(x)=\lambda^3y(x)+f(x)\,(x>0),\quad y(0)=y_0,\quad y'(0)=y_1,\quad y''(0)=y_2\label{eq1.6}
\end{equation}
for $f=0$, due to (v) \eqref{eq1.5}, is
\begin{equation}
y(\lambda,x)=y_0(\lambda,x)=y_0s_0(i\lambda x)+y_1\frac{s_1(i\lambda x)}{i\lambda}+y_2\frac{s_2(i\lambda x)}{(i\lambda)^2}.\label{eq1.7}
\end{equation}
Hence by the method of arbitrary constants one finds solution to the non-homogenous ($f\not=0$) Cauchy problem \eqref{eq1.6}
\begin{equation}
y(\lambda,x)=y_0(\lambda,x)-i\int\limits_0^x\frac{s_2(i\lambda(x-t))}{(i\lambda)^2}f(t)dt,\label{eq1.8}
\end{equation}
here $y_0(\lambda,x)$ is given by \eqref{eq1.7}.

Unit vectors $\{\zeta_k\}_1^3$ \eqref{eq1.3} define three straight lines
\begin{equation}
L_{\zeta_k}\stackrel{\rm def}{=}\{x\zeta_k:x\in\mathbb{R}\}\quad(1\leq k\leq3)\label{eq1.9}
\end{equation}
directions of which are defined by the vectors $\zeta_k$, besides, $\zeta_2L_{\zeta_1}=L_{\zeta_2}$, $\zeta_2L_{\zeta_2}=L_{\zeta_3}$, $\zeta_2L_{\zeta_3}=L_{\zeta_1}$. Denote by $l_{\zeta_k}$ straight half-lines (rays in the direction of unit vectors $\zeta_k$ outgoing from the origin) and by $\widehat{l}_{\zeta_k}$ denote their complements on $L_{\zeta_k}$ (rays coming to the origin),
\begin{equation}
l_{\zeta_k}\stackrel{\rm def}{=}\{x\zeta_k:x\in\mathbb{R}_+\},\quad\widehat{l}_{\zeta_k}\stackrel{\rm def}{=}L_{\zeta_k}\backslash l_{\zeta_k}\quad(1\leq k\leq3).\label{eq1.10}
\end{equation}
Straight lines $\{L_{\zeta_k}\}_1^3$ divide plane $\mathbb{C}$ into six sectors,
\begin{equation}
s_p\stackrel{\rm def}{=}\left\{z\in\mathbb{C}:\frac{2\pi}6(p-1)<\arg z<\frac{2\pi}6p\right\}\quad(1\leq p\leq6).\label{eq1.11}
\end{equation}

\begin{lemma}\label{l1.2}
Zeros of the functions $s_0(z)$, $s_1(z)$, $s_2(z)$ \eqref{eq1.4} lie on the rays $\{\widehat{l}_{\zeta_k}\}$ \eqref{eq1.10} and are given by
\begin{equation}
\{-\zeta_2^lx_0(k)\}_{k=1}^\infty;\quad\{-\zeta_2^lx_1(k)\}_{k=1}^\infty;\quad\{-\zeta_2^lx_2(k)\}_{k=1}^\infty;\label{eq1.12}
\end{equation}
where $l=-1$, $0$, $1$; and $x_0(k)$, $x_1(k)$, $x_2(k)$ are non-negative numbers enumerated in the ascending order. Besides, $x_0(k)$ are positive simple roots of the equation
\begin{equation}
\cos\frac{\sqrt3}2x=-\frac12e^{-3/2x}\quad(x_0(1)>0)\label{eq1.13}
\end{equation}
and $x_1(k)$ and $x_2(k)$ are non-negative simple roots of equations
\begin{equation}
\begin{array}{ccc}
{\displaystyle\cos\left(\frac{\sqrt3}2x-\frac\pi 3\right)=\frac12e^{-3/2x}\quad(x_1(1)=0),}\\
{\displaystyle\cos\left(\frac{\sqrt3}2x+\frac\pi 3\right)=\frac12e^{-3/2x}\quad(x_2(1)=0).}
\end{array}\label{eq1.13'}
\end{equation}

The sequence $\{x_2(k)\}$ interlaces with the sequence $\{x_1(k)\}$ which, in its turn, interlaces with the sequence $\{x_0(k)\}$.
\end{lemma}

Equations \eqref{eq1.13}, \eqref{eq1.13'} readily imply asymptotic behavior of $x_0(k)$, $x_1(k)$, $x_2(k)$ when $k\rightarrow\infty$.
\vspace{5mm}

{\bf 1.2.} The Fourier transform \cite{23, 24}
\begin{equation}
\widetilde{f}(\lambda)=\mathcal{F}(f)(\lambda);\quad\widetilde{f}(\lambda)\stackrel{\rm def}{=}\int\limits_{\mathbb{R}}e^{-i\lambda x}f(x)dx\quad(f\in L^2(\mathbb{R})),\label{eq1.14}
\end{equation}
due to $f=f_++f_-$ ($f_\pm=f\chi_\pm$, $\chi_\pm$ is characteristic function of half-axes $\mathbb{R}_\pm$), is given by
$$\widetilde{f}(\lambda)=\int\limits_0^\infty e^{-i\lambda x}f_+(x)dx+\int\limits_0^\infty e^{i\lambda x}f_-(-x)dx.$$
Hence, $\widetilde{f}(\lambda)$ is equal to the sum of Fourier transforms on rays $\mathbb{R}_\pm$,
\begin{equation}
\widetilde{f}(\lambda)=\widetilde{f}_+(\lambda)+\widetilde{f}_-(\lambda);\quad\widetilde{f}_\pm(\lambda)\stackrel{\rm def}{=}\int\limits_0^\infty e^{-i\lambda\zeta_\pm x}f_\pm(\zeta_\pm x)dx,\label{eq1.15}
\end{equation}
here $\zeta_\pm=\pm1$ and $\{e^{-i\lambda\zeta_\pm x}\}$ are solutions to the equation $-D^2y(x)=\lambda^2y(x)$. Functions $\widetilde{f}_\pm(\lambda)$ \eqref{eq1.15} belong to $L^2(\mathbb{R})$ and $\widetilde{f}_+\perp\widetilde{f}_-$, besides, Parseval's equality holds, $\|f_\pm\|_{L^2(\mathbb{R}_\pm)}=2\pi\|\widetilde{f}_\pm\|_{L^2(\mathbb{R})}$ \cite{23, 24}. Paley -- Wiener theorem \cite{23, 24} implies that $\widetilde{f}_+(\lambda)$ ($\widetilde{f}_-(\lambda)$) is holomorphically extended into the lower $\mathbb{C}_-$ (the upper $\mathbb{C}_+$) half-plane and is of Hardy class $H_-^2$ ($H_+^2$) corresponding to $\mathbb{C}_-$ ($\mathbb{C}_+$). So, $\widetilde{f}(\lambda)$ \eqref{eq1.14} has natural realization as a Fourier transform on the rays \eqref{eq1.15}, here $\{e^{-i\lambda\zeta_\pm x}\}$ are solutions to a second-order equation.

Proceed to the Fourier transform generated by the system of three exponents $\{e^{i\lambda\zeta_kx}\}_1^3$ that are solutions to the third-order equation \eqref{eq1.1}. Consider the bundle
\begin{equation}
l\stackrel{\rm def}{=}\bigcup\limits_kl_{\zeta_k}\label{eq1.16}
\end{equation}
generated by the rays $l_{\zeta_k}$ \eqref{eq1.10} and let $\chi_k=\chi_{l_{\zeta_k}}$ be characteristic functions of the sets $l_{\zeta_k}$ ($1\leq k\leq3$). Denote by $L^2(l)$ the Hilbert space of functions
\begin{equation}
\begin{array}{ccc}
{\displaystyle L^2(l)\stackrel{\rm def}{=}\left\{f=f_1\chi_1+f_2\chi_2+f_3\chi_3,\quad\supp f_k\subseteq l_{\zeta_k}\quad(1\leq k\leq3);\right.}\\
{\displaystyle\left.\|f\|^2=\sum\limits_k\int\limits_0^\infty|f_k(x\zeta_k)|^2dx<\infty\right\}}
\end{array}\label{eq1.17}
\end{equation}
with scalar product
$$\langle f,g\rangle=\sum\limits_k\int\limits_0^\infty f_k(x\zeta_k)\overline{g_k(x\zeta_k)}dx.$$
Match every component $f_k$ ($1\leq k\leq3$) of an element $f\in L^2(l)$ with its Fourier transform
\begin{equation}
\begin{array}{ccc}
{\displaystyle\widetilde{f}_1(\lambda)=(\mathcal{F}_{\zeta_1}f_1)(\lambda)\stackrel{\rm def}{=}\int\limits_0^\infty e^{-i\lambda\zeta_1 x}f_1(x\zeta_1)dx;}\\
{\displaystyle\widetilde{f}_2(\lambda)=(\mathcal{F}_{\zeta_2}f_2)(\lambda)\stackrel{\rm def}{=}\int\limits_0^\infty e^{-i\lambda\zeta_2x}f_2(x\zeta_2)dx;}
\end{array}\label{eq1.18}
\end{equation}
$$\widetilde{f}_3(\lambda)=(\mathcal{F}_{\zeta_3}f_3)(\lambda)\stackrel{\rm def}{=}\int\limits_0^\infty e^{-i\lambda\zeta_3x}f_3(x\zeta_3)dx.$$

Proceed to the description of $\{\widetilde{f}_k(\lambda)\}$ \eqref{eq1.18}.

(A) The function $f_1(\lambda)$ \eqref{eq1.18} is holomorphically extendable into the lower half-plane
\begin{equation}
\mathbb{C}_-(\zeta_1)(=\mathbb{C}_-)\stackrel{\rm def}{=}\{\lambda=\mu+i\nu\in\mathbb{C},\,\nu<0\}\label{eq1.19}
\end{equation}
and due to the Paley -- Wiener theorem \cite{23, 24} the operator $\mathcal{F}_{\zeta_1}$ \eqref{eq1.18} is a unitary isomorphism between $L^2(l_{\zeta_1})$ and the space $H_-^2(\zeta_1)$ ($=H_-^2$), besides, $\|f_1\|_{L^2(f_{\zeta_1})}=2\pi\|\widetilde{f}_1\|^2_{L^2(L\zeta_1)}$. The inverse operator $\mathcal{F}_\zeta^{-1}$ is given by
$$(\mathcal{F}_{\zeta_1}^{-1}g)(x)=\frac1{2\pi}\int\limits_{L_{\zeta_1}}e^{i\lambda\zeta_1x}g(\lambda\zeta_1)d\lambda\quad(g\in H_-^2(\zeta_1))$$
and $\mathcal{F}_{\zeta_1}^{-1}=(\mathcal{F}_{\zeta_1})^*$.

(B) A number $\lambda\zeta_2$ is real if only $\lambda\in L_{\zeta_3}$ \eqref{eq1.9}, i. e., $\lambda=\zeta_3\eta$ ($\eta\in\mathbb{R}$) and $\lambda\zeta_2=\eta\in\mathbb{R}$, therefore
$$\widetilde{f}_2(\lambda)=\int\limits_0^\infty e^{-i\eta x}f_2(x\zeta_2)dx\quad(\lambda=\zeta_3\eta\in L_{\zeta_3})$$
and thus on the straight line $L_{\zeta_3}$ the function $\widetilde{f}_2(\lambda)$ ($\lambda=\zeta_3\eta\in L_{\zeta_3}$) is a usual Fourier transform. Moreover, $\widetilde{f}_2(\lambda)$ has analytic extension into the half-plane
\begin{equation}
\mathbb{C}_-(\zeta_3)\stackrel{\rm def}{=}\{\lambda=\zeta_3\eta,\eta=\mu+i\nu,\,\nu<0\}\label{eq1.20}
\end{equation}
above the straight line $L_{\zeta_3}$, besides, $\lambda\in\mathbb{C}_-(\zeta_3)\Leftrightarrow\eta\in\mathbb{C}_-(\zeta_1)$ \eqref{eq1.19} where $\lambda=\eta\zeta_3$. The straight line
$$L_{\zeta_3}(a)\stackrel{\rm def}{=}\{\lambda=\zeta_3(\eta-ia):\eta,a\in\mathbb{R}\}$$
is parallel to $L_{\zeta_3}$ and for $a>0$ it lies in $\mathbb{C}_-(\zeta_3)$, besides,
$$\widetilde{f}_2(\lambda)=\int\limits_0^\infty e^{-i\eta x}e^{-ax}f_2(\zeta_2x)dx\quad(\forall\lambda\in L_{\zeta_3}(a)).$$
Therefore, integrals of squares of modules of $\widetilde{f}_2(\lambda)$ along every straight line $L_{\zeta_3}(a)$ ($a>0$) are uniformly bounded. Such functions form a Hardy class $H_-^2(\zeta_3)\subset L^2(L_{\zeta_3})$ corresponding to the half-plane $\mathbb{C}_-(\zeta_3)$. So, the operator $\mathcal{F}_{\zeta_2}$ \eqref{eq1.18} is an isomorphism between the spaces $L^2(l_{\zeta_2})$ and $H_-^2(\zeta_3)$, besides, the inverse $\mathcal{F}_{\zeta_2}^{-1}$ is given by
$$(\mathcal{F}_{\zeta_2}^{-1}g)(x)=\frac1{2\pi}\int\limits_{L(\zeta_3)}e^{i\lambda\zeta_2x}g(\lambda\zeta_3)d\lambda\quad(g\in H_-^2(\zeta_3)),$$
and Parseval's identity becomes
$$\langle\mathcal{F}_{\zeta_2}(f_2)(\lambda),\mathcal{F}_{\zeta_2}(\widehat{f}_2)(\lambda\zeta_2)\rangle_{L^2(L_{\zeta_3})}=2\pi\langle f_2(x),\widehat f_2(x)\rangle_{L^2(l_{\zeta_2})}$$
($\forall f_2$, $\widehat f_2\in L^2(l_{\zeta_2})$).

(C) The number $\lambda\zeta_3$ is real if $\lambda\in L_{\zeta_2}$ \eqref{eq1.9}, i. e., $\lambda=\zeta_2\eta$ ($\eta\in\mathbb{R}$) and $\lambda\zeta_3=\eta\in\mathbb{R}$, thus
$$\widetilde{f}_3(\lambda)=\int\limits_0^\infty e^{-i\eta x}f_3(x\zeta_3)dx\quad(\lambda=\zeta_2\eta\in L_{\zeta_2}).$$
Hence it follows that $\widetilde{f}_3(\lambda)$ is holomorphically extendable into the half-plane
\begin{equation}
\mathbb{C}_-(\zeta_2)\stackrel{\rm def}{=}\{\lambda=\zeta_2\eta:\eta=\mu+i\nu,\,\nu<0\}\label{eq1.21}
\end{equation}
above the straight line $L_{\zeta_2}$ and $\lambda\in\mathbb{C}_-(\zeta_2)\Leftrightarrow\eta\in\mathbb{C}_-(\zeta_1)$ \eqref{eq1.19} ($\lambda=\zeta_2\eta$). The function $\widetilde f_3(\lambda)$ belongs to the Hardy space $H_-^2(\zeta_2)\subset L^2(L_{\zeta_2})$ corresponding to the half-plane $\mathbb{C}_-(\zeta_2)$. As in (B), the operator $\mathcal{F}_{\zeta_3}$ \eqref{eq1.18} sets an isomorphism between the spaces $L^2(L_{\zeta_3})$ and $H_-^2(\zeta_2)$ and the inverse to it operator is
$$(\mathcal{F}_{\zeta_3}^{-1}g)(x)=\frac1{2\pi}\int\limits_{L_{\zeta_2}}e^{i\lambda\zeta_3x}g(\lambda\zeta_2)dx\quad(g\in H_-^2(\zeta_2)),$$
besides, Parseval's equality is
$$\langle(\mathcal{F}_{\zeta_3}f_3)(\lambda),(\mathcal{F}_{\zeta_3}\widehat{f}_3)(\lambda\zeta_3)\rangle_{L_{\zeta_2}^2}=2\pi\langle f_3(x),\widehat{f}_3(x)\rangle_{L^2(l_{\zeta_3})}.$$
\begin{picture}(200,200)
\put(0,100){\vector(1,0){200}}
\put(110,100){\vector(1,0){10}}
\put(150,0){\vector(-1,2){100}}
\put(100,100){\vector(-1,2){10}}
\put(100,100){\vector(-1,-2){10}}
\put(150,200){\vector(-1,-2){100}}
\put(30,190){$L_{\zeta_2}$}
\put(65,190){$\mathbb{C}_-(\zeta_2)$}
\put(150,190){$\widetilde{f}_3(\lambda)$}
\put(25,160){$\mathbb{C}_+(\zeta_2)$}
\put(0,105){$\widetilde{f}_1(\lambda)$}
\put(80,110){$\zeta_2$}
\put(90,70){$\zeta_3$}
\put(115,89){$\zeta_1$}
\put(105,103){$0$}
\put(190,105){$L_{\zeta_1}$}
\put(140,105){$\mathbb{C}_+(\zeta_1)$}
\put(140,85){$\mathbb{C}_-(\zeta_1)$}
\put(20,30){$\mathbb{C}_-(\zeta_3)$}
\put(70,30){$\mathbb{C}_+(\zeta_3)$}
\put(140,30){$\widetilde{f}_2(\lambda)$}
\put(60,0){$L_{\zeta_3}$}
\end{picture}

\hspace{20mm}Fig. 1

Define the Hilbert space
\begin{equation}
\begin{array}{ccc}
L^2(L)\stackrel{\rm def}{=}\left\{\widetilde{f}=\sum\widetilde{f}_k\widetilde{\chi}_k:\supp\widetilde{f}_k\subseteq L_{\zeta_k}\,(1\leq k\leq3);\right.\\
{\displaystyle\left.\|f\|^2=\sum\limits_k\int\limits_{L_{\zeta_k}}|\widetilde{f}_k(x\zeta_k)|^2dx<\infty\right\}}
\end{array}\label{eq1.22}
\end{equation}
given on the bundle of straight lines
\begin{equation}
L\stackrel{\rm def}{=}\bigcup\limits_kL_{\zeta_k},\label{eq1.23}
\end{equation}
here $\{L_{\zeta_k}\}$ are from \eqref{eq1.9} and $\widetilde{\chi}_k=\chi_{L_{\zeta_k}}(x)$ are characteristic functions of the sets $L_{\zeta_k}$ ($1\leq k\leq3$).

\begin{theorem}\label{t1.1}
Fourier transform $\{\mathcal{F}_{\zeta_k}\}$ \eqref{eq1.18} is a unitary isomorphism of the space $L^2(l)$ \eqref{eq1.17} onto the subspace $H_-^2(\zeta_1)\oplus H_-^2(\zeta_2)\oplus H_-^2(\zeta_3)$ in $L^2(L)$ \eqref{eq1.22}.
\end{theorem}

For the subspace $L^2(\widehat{l})$, where
\begin{equation}
\widehat{l}=\bigcup\limits_k\widehat{l}_{\zeta_k}\label{eq1.24}
\end{equation}
($\widehat{l}_{\zeta_k}$ are from \eqref{eq1.10}), an analogue of this theorem with natural substitution of Hardy spaces $H_-^2(\zeta_k)\rightarrow H_+^2(\zeta_k)$ ($1\leq k\leq3$) is true. Theorem \ref{t1.1} is a generalization of Paley -- Wiener theorem onto the space of functions given on the bundle of rays $l$ \eqref{eq1.16}.

\begin{remark}\label{r1.2}
Since ${\displaystyle\bigcap\limits_k\mathbb{C}_-(\zeta_k)=\{0\}}$, where $\mathbb{C}_-(\zeta_k)$ are given by \eqref{eq1.19} -- \eqref{eq1.21}, Fourier transforms $\{\widetilde{f}_k(\lambda)\}$ \eqref{eq1.18} exist simultaneously only for $\lambda=0$.
\end{remark}
\vspace{5mm}

{\bf 1.3} Define the space
\begin{equation}
\begin{array}{ccc}
{\displaystyle L^2(l,a)\stackrel{\rm def}{=}\left\{f=\sum f_k\chi_k:\supp f_k\subseteq l_{\zeta_k}\,(1\leq k\leq3);\right.}\\
{\displaystyle\left.\sum\limits_k\int\limits_0^\infty e^{2ax}|f_k(x\zeta_k)|^2dx<\infty\right\}}
\end{array}\label{eq1.25}
\end{equation}
where $a\in\mathbb{R}$. Obviously, $L^2(l,b)\subseteq L^2(l,a)$ for $b\geq a$, in particular, $L^2(l,a)\subseteq L^2(l)$ for all $a\geq0$. Consider the Fourier transform \eqref{eq1.18} of the components of a function $f$ from $L^2(l,a)$ considering that $a\geq0$.

(${\rm A_1}$) The equality
$$\widetilde{f}_1(\lambda)=\int\limits_0^\infty e^{-i\lambda\zeta_1x}f_1(x\zeta_1)=\int\limits_0^\infty e^{-i\zeta_1(\lambda-ia)x}e^{ax}f_1(x_k\zeta_1)dx$$
and $e^{ax}f_1(x\zeta_1)\in L^2(l_{\zeta_1})$ imply that $\widetilde{f}_1(\lambda)$ is holomorphically extendable into the half-plane
\begin{equation}
\mathbb{C}_-(\zeta_1,a)\stackrel{\rm def}{=}\{\lambda=\mu+i\nu\in\mathbb{C}:\nu<a\}\label{eq1.26}
\end{equation}
to which there corresponds the Hardy class $H_-^2(\zeta,a)$ \cite{24}.

(${\rm B_1}$) For $\widetilde{f}_2(\lambda)$ \eqref{eq1.18}, the following equality holds:
$$\widetilde{f}_2(\lambda)=\int\limits_0^\infty e^{-i\lambda\zeta_2x}f_2(x\zeta_2)dx=\int\limits_0^\infty e^{-i\zeta_2(\lambda-\i\zeta_3a)}e^{ax}f_2(x\zeta_2)dx,$$
therefore, due to $e^{ax}f_2(x\zeta_2)\in L^2(l_{\zeta_2})$, hence it follows (see Subseq. (B)) that $\widetilde{f}_2(\lambda)$ is analytically extendable into the half-plane
\begin{equation}
\mathbb{C}_-(\zeta_3,a)\stackrel{\rm def}{=}\{\lambda=\zeta_3\eta:\eta=\mu+i\nu,\nu<a\}\label{eq1.27}
\end{equation}
and thus $\widetilde{f}_2(\lambda)$ belongs to the Hardy space $H_-^2(\zeta_3,a)$ corresponding to the half-plane $\mathbb{C}_-(\zeta_3,a)$ \eqref{eq1.27}.

(${\rm C_1}$) Analogously, $\widetilde{f}_3(\lambda)$ \eqref{eq1.18} for $f\in L^2(l,a)$ has holomorphic extension into the half-plane
\begin{equation}
\mathbb{C}_-(\zeta_2,a)\stackrel{\rm def}{=}\{\lambda=\zeta_2\eta:\eta=\mu+i\nu,\nu<a\}\label{eq1.28}
\end{equation}
and belongs to the Hardy class $H_-^2(\zeta_2,a)$ corresponding to $\mathbb{C}_-(\zeta_2,a)$ \eqref{eq1.28}.
\begin{picture}(200,200)
\put(0,100){\vector(1,0){200}}
\put(0,130){\line(1,0){200}}
\put(150,0){\vector(-1,2){100}}
\put(120,0){\line(-1,2){100}}
\put(180,200){\line(-1,-2){100}}
\put(150,200){\vector(-1,-2){100}}
\put(30,190){$L_{\zeta_2}$}
\qbezier[90](100,0)(100,100)(100,200)
\put(125,190){$\mathbb{C}_-(\zeta_3,a)$}
\put(40,160){$\mathbb{C}_-(\zeta_2,a)$}
\put(190,105){$L_{\zeta_1}$}
\put(142,120){$\mathbb{C}_-(\zeta_1,a)$}
\put(102,130){$a$}
\put(100,40){$-2a$}
\put(87,103){$0$}
\put(110,105){$T_a$}
\put(60,0){$L_{\zeta_3}$}
\end{picture}

\hspace{20mm} Fig. 2

The intersection of the half-planes $\{\mathbb{C}_-(\zeta_k,a)\}$ \eqref{eq1.26} -- \eqref{eq1.28} is an equilateral triangle $T_a$ (see Fig. 2) with side $2\sqrt3a$,
\begin{equation}
T_a\stackrel{\rm def}{=}\bigcap\limits_k\overline{\mathbb{C}_-(\zeta_k,a)}=\{\lambda=\mu+i\nu\in\mathbb{C}:\nu<a,\nu\geq\sqrt3\mu-2a,\nu\geq-\sqrt3\mu-2a\}.\label{eq1.29}
\end{equation}

\begin{theorem}\label{t1.2}
Fourier transforms $\{\mathcal{F}_{\zeta_k}\}$ \eqref{eq1.18} set isomorphism of the space $H_-^2(\zeta_1,a)\oplus H_-^2(\zeta_2,a)\oplus H_-^2(\zeta_3,a)$ into $L^2(L)$ \eqref{eq1.22}. At every point $\lambda$ of the triangle $T_a$ \eqref{eq1.29} all Fourier transforms $\{\widetilde{f}_k(\lambda)\}$ \eqref{eq1.18} exist simultaneously and $\widetilde{f}_k(\lambda)$ ($1\leq k\leq3$) are holomorphic at every inner point $\lambda\in T_a$.
\end{theorem}
\vspace{5mm}

{\bf 1.4} Let $\widetilde{f}_1(\lambda)$ be Fourier transform \eqref{eq1.18} of the first component of $f\in L^2(l,a)$ \eqref{eq1.25} where $a\geq0$; then functions $\widetilde{f}_1(-\lambda)$ define the Hardy space $H_+^2(\zeta_1,-a)$ corresponding to the half-plane
\begin{equation}
\mathbb{C}_+(\zeta_1,-a)\stackrel{\rm def}{=}\{\lambda=\mu+i\nu\in\mathbb{C}:\nu>-a\}.\label{eq1.30}
\end{equation}
In the same way, functions $\widetilde{f}_2(-\lambda)$  ($\widetilde{f}_2(\lambda)$ is Fourier transform \eqref{eq1.18} of the second component $f\in l^2(l,a)$) define the Hardy class $H_+^2(\zeta_3,-a)$ corresponding to the half-plane
\begin{equation}
\mathbb{C}_+(\zeta_3,-a)\stackrel{\rm def}{=}\{\lambda=\zeta_3\eta:\eta=\mu+i\nu,\nu>-a\}.\label{eq1.31}
\end{equation}
Finally, $\widetilde{f}_3(-\lambda)$ ($\widetilde{f}_3(\lambda)$ is Fourier transform \eqref{eq1.18} of the third component of $f\in L^2(l,a)$) generate the Hardy space $H_+^2(\zeta_2,-a)$ corresponding to the half-plane
\begin{equation}
\mathbb{C}_+(\zeta_2,a)\stackrel{\rm def}{=}\{\lambda=\zeta_2\eta:\eta=\mu+i\nu,\nu>-a\}.\label{eq1.32}
\end{equation}
\begin{picture}(200,200)
\put(0,100){\vector(1,0){200}}
\put(0,70){\line(1,0){200}}
\put(150,0){\vector(-1,2){100}}
\put(180,0){\line(-1,2){100}}
\put(120,200){\line(-1,-2){100}}
\put(150,200){\vector(-1,-2){100}}
\put(30,190){$L_{\zeta_2}$}
\qbezier[90](100,0)(100,100)(100,200)
\put(125,190){$\mathbb{C}_+(\zeta_3,-a)$}
\put(27,175){$\mathbb{C}_+(\zeta_2,-a)$}
\put(190,105){$L_{\zeta_1}$}
\put(142,80){$\mathbb{C}_+(\zeta_1,-a)$}
\put(102,160){$2a$}
\put(100,70){$-a$}
\put(87,103){$0$}
\put(75,80){$T_a^*$}
\put(60,0){$L_{\zeta_3}$}
\end{picture}

\hspace{20mm} Fig. 3

Intersection of half-planes $\{\mathbb{C}_+(\zeta_k,-a)\}$ \eqref{eq1.30} -- \eqref{eq1.32} gives an equilateral triangle $T_a^*$ (see Fig. 3) with side $2\sqrt3a$,
\begin{equation}
T_a^*\stackrel{\rm def}{=}\bigcap\limits_k\mathbb{C}_+(\zeta_k,-a)=\{\lambda=\mu+i\nu\in\mathbb{C}:\nu\geq-a;\nu\leq\sqrt3\mu-2a,\nu\leq\sqrt3\mu+2a\}.\label{eq1.33}
\end{equation}
Triangle $T_a^*$ is complexly adjoint to $T_a$ \eqref{eq1.29},
$$\lambda\in T_a\Leftrightarrow\overline{\lambda}\in T_a^*.$$

\begin{theorem}\label{t1.3}
For every function $f\in L^2(l,a)$ \eqref{eq1.25} ($a>0$) for all $\lambda\in T_a\cap T_a^*$ where $T_a$ and $T_a^*$ are given by \eqref{eq1.29} and \eqref{eq1.33}, there exist Fourier transforms $\{\widetilde{f}_k(\lambda)\}$ \eqref{eq1.18} and $\{\widetilde{f}_k(-\lambda)\}$ which are analytic functions at all inner points of the set $T_a\cap T_a^*$.
\end{theorem}

Notice that intersection $T_a\cap T_a^*$ forms regular hexagon with side $a2/\sqrt3$ (see Figs. 2, 3).

\begin{remark}\label{r1.3}
Instead of $L^2(l,a)$ \eqref{eq1.25}, one can consider the space
\begin{equation}
L^2(l,a_k)\stackrel{\rm def}{=}\left\{f=\sum f_k\chi_k: \supp f_k\subseteq l_{\zeta_k};\sum\limits_k\int\limits_0^\infty e^{2a_kx}|f_k(x\zeta_k)|^2dx<\infty\right\}\label{eq1.34}
\end{equation}
where $a_k\in\mathbb{R}$ ($1\leq k\leq3$). For $a_k\geq0$ ($\forall k$) analogues of Theorems \ref{t1.2}, \ref{t1.3} are true but in this case the triangle $T_a$ \eqref{eq1.29} is no longer equilateral.
\end{remark}
\vspace{5mm}

{\bf 1.5}. Consider the limit case ($a=\infty$). Let
\begin{equation}
L^2(l,\infty)\stackrel{\rm def}{=}\bigcap\limits_{a\geq0}L^2(l,a)=\left\{f\in L^2(l):\sum\limits_k\int\limits_0^\infty e^{2ax}|f(x\zeta_k)|^2dx<\infty,\forall a\geq0\right\},\label{eq1.35}
\end{equation}
besides, $L^2(l,\infty)\subseteq L^2(l)$ \eqref{eq1.17} and topology in $L^2(l,\infty)$ is set by the metric of space $L^2(l)$.

(${\rm A_1}$) Functions $\widetilde{f}_1(\lambda)$ \eqref{eq1.18}, where $f_1$ is the first component of $f\in L^2(l,\infty)$, form the {\bf space of entire functions}
\begin{equation}
HE_-^2(\zeta_1)=H_-^2(\zeta_1,\infty)\stackrel{\rm def}{=}\bigcap\limits_{a\geq0}H_-^2(\zeta_1,a).\label{eq1.36}
\end{equation}
Space $HE_-^2(\zeta_1)$ consists of the functions from $H_-^2(\zeta_1)$ that are holomorphically extendable into the half-plane
$$\mathbb{C}_+(\zeta_1)\stackrel{\rm def}{=}\{\lambda=\mu+i\nu\in\mathbb{C}:\nu>0\}$$
and squares of modules of which are summable on every horizontal straight line in $\mathbb{C}_+(\zeta_1)$.

(${\rm B_1}$) Second components $\widetilde{f}_2(\lambda)$ \eqref{eq1.18} ($f\in L^2(l,\infty)$) give the {\bf space of entire functions}
\begin{equation}
HE_-^2(\zeta_3)=H_-^2(\zeta_3,\infty)\stackrel{\rm def}{=}\bigcap\limits_{a\geq0}H_-^2(\zeta_3,a),\label{eq1.37}
\end{equation}
topology of which is induced by $H_-^2(\zeta_3)$. The space $HE_-^2(\zeta_3)$ is formed by such functions from $H_-^2(\zeta_3)$ that are analytically extendable into the half-plane
$$\mathbb{C}_+(\zeta_3)\stackrel{\rm def}{=}\{\lambda=\zeta_3\eta:\eta=\mu+i\nu\in\mathbb{C},\nu>0\},$$
and squares of modules of which are integrable on every straight line from $\mathbb{C}_+(\zeta_3)$ that is parallel to $L_{\zeta_3}$.

(${\rm C_1}$) Finally, the functions $\widetilde{f}_3(\lambda)$ \eqref{eq1.18} ($f\in L^2(l,\infty)$) give the {\bf space of entire functions}
\begin{equation}
HE_-^2(\zeta_2)=H_-^2(\zeta_2,\infty)\stackrel{\rm def}{=}\bigcap\limits_{a\geq0}H_-^2(\zeta_2,a)\label{eq1.38}
\end{equation}
holomorphically extendable into the half-plane
$$C_+(\zeta_2)\stackrel{\rm def}{=}\{\lambda=\zeta_2\eta:\eta=\mu+i\nu\in\mathbb{C},\nu>0\}$$
and squares of their modules are summable on every straight line in $\mathbb{C}_+(\zeta_2)$ parallel to $L_{\zeta_2}$.

\begin{theorem}\label{t1.4}
Fourier transform $\{\mathcal{F}_{\zeta_k}\}$ \eqref{eq1.18} is an isomorphism of the space $L^2(l,\infty)$ \eqref{eq1.35} onto the space of entire functions $HE_-^2(\zeta_1)\oplus HE_-^2(\zeta_2)\oplus HE_-^2(\zeta_3)$ ($\subset L^2(L)$ \eqref{eq1.22}) where $\{HE_-^2(\zeta_k)\}$ are given by \eqref{eq1.36} -- \eqref{eq1.38}.
\end{theorem}

\begin{definition}\label{d1.1}
For a function $f\in L^2(l,a)$ \eqref{eq1.25} for each $\lambda\in T_a$ \eqref{eq1.29}, expression
\begin{equation}
\widetilde{f}(\lambda)=(\mathcal{F}(f))(\lambda);\quad(\mathcal{F}f)(\lambda)\stackrel{\rm def}{=}\sum\limits_k\widetilde{f}_k(\lambda)\label{eq1.39}
\end{equation}
is said to be the Fourier transform of $f$ where $\widetilde{f}_k(\lambda)=(\mathcal{F}_{\zeta_k}(f_k))(\lambda)$ \eqref{eq1.18} ($1\leq k\leq3$).
\end{definition}
\vspace{5mm}

{\bf 1.6.} Symmetry $x\rightarrow-x$ of straight line $\mathbb{R}$ ($\pi$ rotation) sets an automorphism $(Jf)(x)=f(-x)$ of the space $L^2(\mathbb{R})$ ($f\in L^2(\mathbb{R}$) which is an involution, $J^2=I$. The space $L^2(\mathbb{R})$ splits into the orthogonal sum, $L^2(\mathbb{R})=E_+\oplus E_-$, besides, $E_+$ ($E_-$) is formed by even (odd) functions corresponding to the eigenvalue $\lambda_+=1$ ($\lambda_-=-1$) of the operator $J$. Fourier transform \eqref{eq1.14} of the subspaces $E_+$ and $E_-$ leads \eqref{eq1.15} to cosine and sine Fourier transforms \cite{23, 24}.

Give an analogue of this property for the space $L^2(l,a)$ \eqref{eq1.25} ($a\geq0$). Functions ${\displaystyle f=\sum f_k\chi_k\in L^2(l,a)}$ are given on a bundle of rays $l$ \eqref{eq1.16}, therefore it is natural to define an automorphism $J$ as $2\pi/3$ rotation,
\begin{equation}
J(f_1(x\zeta_1)\chi_1+f_2(x\zeta_2)\chi_2+f_3(x\zeta_3)\chi_3)=f_3(x\zeta_1)\chi_1+f_1(x\zeta_2)\chi_2+f(x\zeta_3)\chi_3,\label{eq1.40}
\end{equation}
here $\chi_k$ are characteristic functions of rays $l_{\zeta_k}$ \eqref{eq1.10} ($1\leq k\leq3$). Operator $J$ is an isomorphism and $J^3=I$, it has three eigenvalues $\zeta_1$, $\zeta_2$, $\zeta_3$ \eqref{eq1.2} to which there correspond three proper subspaces,
$$E_{\zeta_1}=\{\Phi=\varphi(x\zeta_1)\chi_1+\varphi(x\zeta_2)\chi_2+\varphi(x\zeta_3)\chi_3\in L^2(l,a)\};$$
\begin{equation}
E_{\zeta_2}=\{\Psi=\psi(x\zeta_1)\chi_1+\zeta_3\psi(x\zeta_2)\chi_2+\zeta_2\psi(x\zeta_3)\chi_3\in L^2(l,a)\};\label{eq1.41}
\end{equation}
$$E_{\zeta_3}=\{H=h(x\zeta_1)\chi_1+\zeta_2h(x\zeta_2)\chi_2+\zeta_3h(x\zeta_3)\chi_3\in L^2(l,a)\}$$
besides,
$$J\Phi=\zeta_\Phi,\quad J\Psi=\zeta_2\Psi,\quad JH=\zeta_3H,$$
and $E_{\zeta_k}\perp E_{\zeta_k}$ ($k\not=s$). The decomposition $L^2(l,a)=E_{\zeta_1}\oplus E_{\zeta_2}\oplus E_{\zeta_3}$ is true, so, for an arbitrary function ${\displaystyle f=\sum f_k(x\zeta_k)\chi_k\in L^2(l,a)}$,
$$f=\Phi(f)+\Psi(f)+H(f),$$
here $\Phi(f)\in E_{\zeta_1}$, $\Psi(f)\in E_{\zeta_2}$, $H(f)\in E_{\zeta_3}$ \eqref{eq1.41}, besides,
$$(\varphi(f))(x\zeta_1)=\frac13(f_1(x\zeta_1)+f_2(x\zeta_1)+f_3(x\zeta_1));$$
\begin{equation}
(\psi f)(x\zeta_1)=\frac13(f(x\zeta_1)+\zeta_2f_2(x\zeta_1)+\zeta_3f_3(x\zeta_1));\label{eq1.42}
\end{equation}
$$h(f)(x\zeta_1)=\frac13(f_1(x\zeta_1)+\zeta_3f_2(x\zeta_1)+\zeta_2f_3(x\zeta_1)).$$

{\bf Case 1.} Let $f_1(x\zeta_1)=g(x)$, $f_2(x\zeta_1)=g(x)=f_3(x\zeta_1)$, then $(\varphi(f))(x\zeta_1)=g(x)$, $(\psi f)(x\zeta_1)=(H(f))(x\zeta_1)=0$ and due to \eqref{eq1.18}, \eqref{eq1.39}
\begin{equation}
\begin{array}{ccc}
{\displaystyle(\mathcal{F}(f))(\lambda)=\int\limits_0^\infty e^{-i\lambda\zeta_1x}g(x)dx+\int\limits_0^\infty e^{-i\lambda\zeta_2x}g(x)dx+\int\limits_0^\infty e^{-i\lambda\zeta_3x}g(x)dx}\\
{\displaystyle=3\int\limits_0^\infty s_0(-i\lambda x)g(x)dx}
\end{array}\label{eq1.43}
\end{equation}
where $s_0(z)$ is given by \eqref{eq1.4} since $g(x)=g(x\zeta_k)$ ($\forall k$).

{\bf Case 2.} Let $f_1(x\zeta_1)=g(x)$ and $f_2(x\zeta_1)=\zeta_3g(x)$, $f_3(x\zeta_1)=\zeta_2g(x)$, $f_3(x\zeta_1)=\zeta_2g(x)$, in this case $(\psi(f))(x\zeta_1)=g(x)$ and $(\varphi(f))(x\zeta_1)=(H(f))(x\zeta_1)=0$, and thus
\begin{equation}
\begin{array}{ccc}
{\displaystyle(\mathcal{F}(f))(\lambda)=\int\limits_0^\infty e^{-i\lambda\zeta_1x}g(x)dx+\zeta_3\int\limits_0^\infty e^{-i\lambda\zeta_2x}g(x)dx+\zeta_2\int\limits_0^\infty e^{-i\lambda\zeta_3x}g(x)dx}\\
{\displaystyle=3\int\limits_0^\infty s_1(-i\lambda x)g(x)dx,}
\end{array}\label{eq1.44}
\end{equation}
here $s_1(z)$ is given by \eqref{eq1.4}.

{\bf Case 3.} If $f_1(x,\zeta_1)=g(x)$ and $f_2(x\zeta_1)=\zeta_2g(x)$, $f_3(x\zeta_1)=\zeta_3g(x)$, then $(Hf)(x\zeta_1)=g(x)$ and $(\varphi(f))(x\zeta_1)=(\psi(f))(x\zeta_1)=0$, therefore
\begin{equation}
(\mathcal{F}(f))(\lambda)=3\int\limits_0^\infty s_2(-i\lambda x)g(x)dx\label{eq1.45}
\end{equation}
where $s_2(z)$ is given by \eqref{eq1.4}.

So, Fourier transform $\mathcal{F}$ \eqref{eq1.18}, \eqref{eq1.39} of the space $E_{\zeta_1}$ \eqref{eq1.41} is expressed via the function $s_0(-i\lambda x)$, Fourier transform of the subspace $E_{\zeta_2}$, correspondingly, is expressed via $s_1(-i\lambda x)$, and, finally, Fourier transform of $E_{\zeta_3}$ is expressed via $s_2(-i\lambda x)$. Equations \eqref{eq1.43} -- \eqref{eq1.45} are analogues of cosine and sine Fourier transforms for the exponents $\{e^{i\lambda\zeta_kx}\}_1^3$.

\section{Jost solutions}\label{s2}

{\bf 2.1.} Denote by $\mathcal{H}$ the Hilbert space
\begin{equation}
\mathcal{H}=L^2(\mathbb{R}_-)\oplus L^2(\mathbb{R}_+)\stackrel{\rm def}{=}\{y=(v,u):v\in L^2(\mathbb{R}_-),\,u\in L^2(\mathbb{R}_+)\}.\label{eq2.1}
\end{equation}
Consider the linear operator $\mathcal{L}_\alpha$ in $\mathcal{H}$,
\begin{equation}
\mathcal{L}_\alpha y\stackrel{\rm def}{=}(-iDv,iD^3u+\alpha\langle u,q\rangle q),\label{eq2.2}
\end{equation}
here $y=(v,u)\in\mathcal{H}$; ${\displaystyle D=\frac d{dx}}$; $\alpha\in\mathbb{R}$; $q\in L^2(\mathbb{R}_+)$ and for some $a\geq0$ meets the condition
\begin{equation}
\int\limits_{\mathbb{R}_+}|q(x)|^2e^{2ax}dx<\infty.\label{eq2.3}
\end{equation}
Domain $\mathfrak{D}(\mathcal{L}_\alpha)$ of the operator $\mathcal{L}_\alpha$ is
\begin{equation}
\mathfrak{D}(\mathcal{L}_\alpha)\stackrel{\rm def}{=}\{y=(v,u)\in\mathcal{H}:v\in W_2^1(\mathbb{R}_-),u\in W_2^3(\mathbb{R}_+);u(0)=0,u'(0)=v(0)\}.\label{eq2.4}
\end{equation}
It is easy to see that the operator $\mathcal{L}_\alpha$ \eqref{eq2.2}, \eqref{eq2.4} is self-adjoint.

\begin{remark}\label{r2.1}
Obviously, inequality \eqref{eq2.3} is true for all $b$ such that $0\leq b<a$. Moreover, $e^{cx}q(x)\in L^1(\mathbb{R}_+)$ if $c$ lies in the interval $0\leq c<a$. Really,
$$\int\limits_{\mathbb{R}_+}e^{cx}|q(x)|dx=\int\limits_{\mathbb{R}_+}e^{ax}|q(x)|\cdot e^{(c-a)x}dx<\infty$$
due to the Cauchy -- Schwarz  inequality since $e^{ax}|q(x)|\in L^2(\mathbb{R}_+)$ \eqref{eq2.3} and $e^{(c-a)x}\in L^2(\mathbb{R}_+)$ ($c<a$).
\end{remark}

\begin{remark}\label{r2.2}
Symmetric operator
$$iD^3+\alpha\langle .,q\rangle q$$
acting in $L^2(\mathbb{R}_+)$, domain of which is $\mathfrak{D}=\{u\in W_2^3(\mathbb{R}_+):u(0)=0,u'(0)=0\}$ does not have self-adjoint extensions in $L^2(\mathbb{R}_+)$ \cite{16, 17, 18}. Such extensions are constructed with departure from the space $L^2(\mathbb{R}_+)$. The operator $\mathcal{L}_\alpha$ \eqref{eq2.2}, \eqref{eq2.4} in $\mathcal{H}$ \eqref{eq2.1} is one of such extensions.
\end{remark}

Consider a system of equations in $\mathcal{H}$,
\begin{equation}
\left\{
\begin{array}{lll}
-iDv=\lambda^3v\quad(x\in\mathbb{R}_-)\\
iD^3u+\alpha\langle u,v\rangle v=\lambda^3u\quad(x\in\mathbb{R}_+),
\end{array}\right.\label{eq2.5}
\end{equation}
then the first equation has the solution $e^{i\lambda^3x}$ ($x\in\mathbb{R}_-$), and the second as $x\rightarrow\infty$, due to \eqref{eq2.3}, transforms into equation \eqref{eq1.1} which has three linearly independent solutions $\{e^{i\lambda\zeta_kx}\}_1^3$. Hence, asymptotic behavior of the solution $u(x)$ as $x\rightarrow\infty$ is described by a linear combination of the exponents $\{e^{i\lambda\zeta_kx}\}_1^3$.

Denote by $\psi_k(\lambda,x)$ a {\bf Jost solution} \cite{1, 2, 3, 4, 5, 21, 25} to the equation
\begin{equation}
iD^3u+\alpha\langle u,q\rangle q=\lambda^3u\quad(x\in\mathbb{R}_+)\label{eq2.6}
\end{equation}
satisfying the boundary condition
\begin{equation}
\lim\limits_{x\rightarrow\infty}u(\lambda,x)e^{-i\lambda\zeta_kx}=b_k(\lambda)\quad(1\leq k\leq3),\label{eq2.7}
\end{equation}
here $b_k(\lambda)$ is some function of $\lambda$.

Consider the function
\begin{equation}
\psi_1(\lambda,x)=b_1(\lambda)e^{i\lambda\zeta_1x}-i\int\limits_x^\infty\frac{s_2(i\lambda(x-t))}{(i\lambda)^2}\alpha\langle\psi_1,q\rangle q(t)dt\label{eq2.8}
\end{equation}
where $s_2(z)$ is given by \eqref{eq1.4}. Integral in \eqref{eq2.8} converges uniformly by $x$ for all $\lambda\in T_a$ \eqref{eq1.29} due to \eqref{eq2.3}. The function $\psi_1(\lambda,x)$ is a solution to equation \eqref{eq2.6} (see \eqref{eq1.8}). Proceed to the boundary condition \eqref{eq2.7} for $\psi_1(\lambda,x)$ \eqref{eq2.8}, to do this, consider
$$\psi_1(\lambda,x)e^{-i\lambda\zeta_1x}=b_1(\lambda)+\frac{\alpha i}{3\lambda^2}\langle \psi_1,q\rangle\int\limits_x^\infty e^{-i\lambda\zeta_1x}\left[e^{i\lambda\zeta_1(x-t)}+\zeta_2e^{i\lambda\zeta_2(x-t)}+\zeta_3e^{i\lambda\zeta_3(x-t)}\right]$$
$$\times q(t)dt=b_1+\frac{\alpha i}{3\lambda^2}\langle\psi_1,q\rangle\int\limits_x^\infty\left[e^{-i\lambda\zeta_1t}+\zeta_2e^{\sqrt3\zeta_3\lambda x}e^{-i\lambda\zeta_2t}+\zeta_3e^{-\sqrt3\zeta_2\lambda x}e^{-i\lambda\zeta_3t}\right]q(t)dt.$$
Since
\begin{equation}
\begin{array}{ccc}
{\displaystyle\left|\int\limits_x^\infty\left[e^{-i\lambda\zeta_1t}+\zeta_2e^{\sqrt3\zeta_3\lambda x}e^{-i\lambda\zeta_2t}+e^{-\sqrt3\zeta_2\lambda x}e^{-i\lambda\zeta_3t}\right]q(t)dt\right|\leq\int\limits_x^\infty\left(1+2e^{\sqrt3|\lambda|x}\right)}\\
{\displaystyle \times e^{|\lambda|t}|q(t)|dt\leq\int\limits_x^\infty\left(e^{-i\sqrt3|\lambda|t}+1\right)e^{(1+\sqrt3)|\lambda|t}|q(t)|dt\leq3\int\limits_x^\infty e^{(1+\sqrt3)|\lambda|t}|q(t)|dt}
\end{array}\label{eq2.9}
\end{equation}
and $e^{(1+\sqrt3)|\lambda|x}q(x)\in L^1(\mathbb{R}_+)$ for $(1+\sqrt3)|\lambda|<a$ (Remark \ref{r2.1}), then this integral tends to zero as $x\rightarrow\infty$ if only $|\lambda|<a(1+\sqrt3)^{-1}$. The disc
\begin{equation}
\mathbb{D}_{a/3}\stackrel{\rm def}{=}\{\lambda\in\mathbb{C}:|\lambda|<a/3\}\label{eq2.10}
\end{equation}
lies inside the disc $|\lambda|<a(1+\sqrt3)^{-1}$ which, in its turn, belongs to the triangle $T_a$ \eqref{eq1.29}. Thus, for all $\lambda\in\mathbb{D}_{a/3}$ integral \eqref{eq2.8} converges and the boundary condition \eqref{eq2.7} is satisfied.

Upon scalar multiplying equation \eqref{eq2.8} by $q$, one obtains
\begin{equation}
\langle\psi,q\rangle(1+\alpha im_d(\lambda))=b_1(\lambda)\widetilde{q}_1^*(\lambda),\label{eq2.11}
\end{equation}
here
\begin{equation}
\begin{array}{lll}
\widetilde{q}_k(\lambda)\stackrel{\rm def}{=}\langle e^{-i\lambda\zeta_kx},\overline{q}(x)\rangle;\quad\widetilde{q}_k^*(\lambda)\stackrel{\rm def}{=}\overline{\widetilde{q}_k(\overline{\lambda})}\quad(1\leq k\leq3);\\
{\displaystyle m_{s_2}(\lambda)\stackrel{\rm def}{=}\left\langle\int\limits_x^\infty\frac{s_2(i\lambda(x-t))}{(i\lambda)^2}q(t)dt,q(x)\right\rangle}
\end{array}\label{eq2.12}
\end{equation}
and due to the definition of $s_2(z)$ \eqref{eq1.4}
\begin{equation}
m_{s_2}(\lambda)=\frac1{3(i\lambda)^2}\sum\zeta_km_k(\lambda),\label{eq2.13}
\end{equation}
besides,
\begin{equation}
m_k(\lambda)\stackrel{\rm def}{=}\left\langle\int\limits_x^\infty e^{i\lambda\zeta_k(x-t)}q(t)dt,q(x)\right\rangle=\int\limits_{\mathbb{R}_+}e^{-i\lambda\zeta_ks}g_q(s)ds\quad(1\leq k\leq3);\label{eq2.14}
\end{equation}
\begin{equation}
g_q(s)\stackrel{\rm def}{=}\int\limits_{\mathbb{R}_+}q(x+s)\overline{q(x)}dx.\label{eq2.15}
\end{equation}

\begin{remark}\label{r2.3}
Since $\mathbb{D}_{a/3}\subset T_a\cap T_a^*$, then the functions $\{q_k(\lambda)\}$ and $\{q_k^*(\lambda)\}$ are defined for all $\lambda\in\mathbb{D}_{a/3}$ (Theorem \ref{t1.2}). The inequality
$$e^{as}|g_q(s)|\leq\int\limits_0^\infty|q(x+s)|e^{a(s+x)}\cdot|\overline{q(x)}|dx$$
and inclusions $q(x)e^{ax}\in L^2(\mathbb{R}_+)$ (see \eqref{eq2.3}), $q\in L^1(\mathbb{R}_+)$ (Remark \ref{r2.1}) imply \cite{23, 24} that $a^sg_q(s)\in L^2(\mathbb{R}_+)$. So, class of functions $q\in L^2(\mathbb{R}_+)$ for which \eqref{eq2.3} takes place is closed with respect to the convolution operation \eqref{eq2.14}, therefore the functions $\{m_{s_2}(\lambda)\}$ \eqref{eq2.14}, and so $m_{s_2}(\lambda)$ \eqref{eq2.12} also, are correctly defined for all $\lambda\mathbb{D}_{a/3}$.
\end{remark}

By $LD(k,a)$ denote union of straight line $L_{\zeta_k}$ \eqref{eq1.9} and disk $\mathbb{D}_{a/3}$ \eqref{eq2.10},
\begin{equation}
LD(k,a)\stackrel{\rm def}{=}L_{\zeta_k}\cup\mathbb{D}_{a/3}\quad(1\leq k\leq3),\label{eq2.16}
\end{equation}
besides, it is obvious that $\zeta_2LD(k,a)=LD((k+1)\mod3,a)$.

\begin{lemma}\label{l2.1}
For the functions $m_k(\lambda)$ \eqref{eq2.14}, the following equalities hold:

$$({\rm {i}})\,m_1(\lambda)+m_1^*(\lambda)=\widetilde{q}_1(\lambda)\widetilde{q}_1^*(\lambda)\quad(\forall\lambda\in LD(1,a));$$
\begin{equation}
({\rm{ii}})\,m_2(\lambda)+m_2^*(\lambda)=\widetilde{q}_2(\lambda)\widetilde{q}_3^*(\lambda)\quad(\forall\lambda\in LD(2,a));\label{eq2.17}
\end{equation}
$$({\rm iii})\,m_3(\lambda)+m_3^*(\lambda)=\widetilde{q}_3(\lambda)\widetilde{q}_2^*(\lambda)\quad(\forall\lambda\in LD(3,a));$$
here $\widetilde{q}_k(\lambda)$ are given by \eqref{eq2.12} and
\begin{equation}
m_1(\lambda\zeta_2)=m_2(\lambda);\quad m_2(\lambda\zeta_2)=m_3(\lambda);\quad m_3(\lambda\zeta_2)=m_1(\lambda).\label{eq2.18}
\end{equation}
\end{lemma}

Proof of the lemma follows from the formulas
$$m_k(\lambda)=\int\limits_0^\infty\int\limits_x^\infty e^{-i\lambda\zeta_kt}q(t)dte^{i\lambda\zeta_kx}\overline{q(x)}dx\quad(1\leq k\leq3)$$
upon integration by parts.

\begin{corollary}\label{c2.1}
Functions $m_{s_2}(\lambda)$ \eqref{eq2.12} and
\begin{equation}
m_{s_0}(\lambda)\stackrel{\rm def}{=}\left\langle\int\limits_x^\infty s_0(i\lambda(x-t)q(t)dt,q(x)\right\rangle;\quad m_{s_1}(\lambda)\stackrel{\rm def}{=}\left\langle\int\limits_x^\infty\frac{s_1(i\lambda(x-t))}{i\lambda}q(t)dt,q(x)\right\rangle\label{eq2.19}
\end{equation}
for all $\lambda\in\mathbb{D}_{a/3}$ \eqref{eq2.10} satisfy the relations
$$m_{s_2}(\lambda)+m_{s_2}^*(\lambda)=\frac1{3(i\lambda)^2}\{\zeta_1\widetilde{q}_1(\lambda)\widetilde{q}_1^*(\lambda)+\zeta_2\widetilde{q}_2(\lambda)\widetilde{q}_3^*(\lambda)+\zeta_3\widetilde{q}_3(\lambda)\widetilde{q}_2^3(\lambda)\};$$
\begin{equation}
m_{s_1}(\lambda)+m_{s_1}^*(\lambda)=\frac1{3(i\lambda)}\{\zeta_1\widetilde{q}_1(\lambda)\widetilde{q}_1^*(\lambda)+\zeta_3\widetilde{q}_2(\lambda)\widetilde{q}_3^*(\lambda)+\zeta_2\widetilde{q}_3(\lambda)\widetilde{q}_2^*(\lambda)\};
\label{eq2.20}
\end{equation}
$$m_{s_0}(\lambda)+m_{s_0}^*(\lambda)=\frac13\{\widetilde{q}_1(\lambda)\widetilde{q}_1^*(\lambda)+\widetilde{q}_2(\lambda)\widetilde{q}_3^*(\lambda)+\widetilde{q}_3(\lambda)\widetilde{q}_2^*(\lambda)\},$$
besides,
\begin{equation}
m_{s_p}(\lambda\zeta_2)=m_{s_p}\quad(p=0,1,2).\label{eq2.21}
\end{equation}
\end{corollary}

Supposing that $b_1(\lambda)=1+\alpha im_{s_2}(\lambda)$ in \eqref{eq2.11}, one finds $\langle\psi_1,q\rangle=q_1^*(\lambda)$, therefore
\begin{equation}
\psi_1(\lambda,x)=b(\lambda)e^{i\lambda\zeta_1x}-\alpha i\widetilde{q}_1^*(\lambda)\int\limits_x^\infty\frac{s_2(i\lambda(x-t))}{(i\lambda)^2}q(t)dt\quad(\lambda\in\mathbb{D}_{a/3})\label{eq2.22}
\end{equation}
where
\begin{equation}
b(\lambda)=1+\alpha im_{s_2}(\lambda)\label{eq2.23}
\end{equation}
and $m_{s_2}(\lambda)$ is given by \eqref{eq2.13}. Similarly, it is proved that the rest of Jost solutions are
\begin{equation}
\begin{array}{ccc}
{\displaystyle\psi_2(\lambda,x)=b(\lambda)e^{i\lambda\zeta_2x}-\alpha i\widetilde{q}_3^*(\lambda)\int\limits_x^\infty\frac{s_2(i\lambda(x-t))}{(i\lambda)^2}q(t)dt\quad(\lambda\in\mathbb{D}_{a/3});}\\
{\displaystyle\psi_3(\lambda,x)=b(\lambda)e^{i\lambda\zeta_3x}-\alpha i\widetilde{q}_2^*(\lambda)\int\limits_x^\infty\frac{s_2(i\lambda(x-t))}{(i\lambda)^2}q(t)dt\quad(\lambda\in\mathbb{D}_{a/3}).}
\end{array}\label{eq2.22}
\end{equation}

\begin{theorem}\label{t2.1}
If $q$ from $L^2(\mathbb{R}_+)$ satisfies condition \eqref{eq2.3}, then Jost solutions $\{\psi_k(\lambda,x)\}$ of equation \eqref{eq2.6} exist for all $\lambda\in\mathbb{D}_{a/3}$ \eqref{eq2.10} and are given by \eqref{eq2.22}. The functions $\{\psi_k(\lambda,x)\}$ are analytic by $\lambda$ for $\lambda\in\mathbb{D}_{a/3}$ \eqref{eq2.10} and
\begin{equation}
\psi_1(\lambda\zeta_2,x)=\psi_2(\lambda,x);\quad\psi_2(\lambda\zeta_2,x)=\psi_3(\lambda\zeta_2,x)=\psi_1(\lambda,x).\label{eq2.24}
\end{equation}
\end{theorem}

\begin{lemma}\label{l2.2}
Jost solutions $\{\psi_k(\lambda,x)\}_1^3$ \eqref{eq2.22} are
\begin{equation}
\psi_k(\lambda,x)=e^{i\lambda\zeta_kx}\{b(\lambda)+\varphi_k(\lambda,x)\}\quad(\lambda\in\mathbb{D}_{a/3}),\label{eq2.25}
\end{equation}
besides, $\varphi_k(\lambda,x)\in L^2(\mathbb{R}_+)$ for all $\lambda$ from $\mathbb{D}_{a/3}$ ($1\leq k\leq3$).
\end{lemma}

P r o o f. Give proof for $\psi_1(\lambda,x)$ (for $\psi_k(\lambda,x)$, $k\not=1$ considerations are similar). Since
$$\psi_1(\lambda,x)e^{-i\lambda\zeta_1x}-b(\lambda)=\varphi_1(\lambda,x)=\frac{\alpha i}{\lambda^2}\widetilde{q}_1^*(\lambda)\int\limits_x^\infty e^{-i\lambda\zeta_1x}s_2(i\lambda(x-t))q(t)dt,$$
then one has to show that the integral as a function of $x$ belongs to $L^2(\mathbb{R}_+)$ when $\lambda\in\mathbb{D}_{a/3}$. To do this, according to estimate \eqref{eq2.9}, it is sufficient to ascertain that the function (of $x$)
$$c(x)=\int\limits_x^\infty e^{(1+\sqrt3)|\lambda|t}|q(t)|dt$$
is square summable. Use the following Hardy theorem \cite{23}.

\begin{theorem}[Hardy]\cite{23}
If $f\in L^2(\mathbb{R}_+)$, then the functions
$$\frac1x\int\limits_0^xf(t)dt;\quad\int\limits_x^\infty\frac1tf(t)dt$$
also belong to $L^2(\mathbb{R}_+)$.
\end{theorem}

The Hardy theorem implies that in order that the function
\begin{equation}
\varepsilon c(x)=\int\limits_x^\infty\frac 1te^{(1+\sqrt3)|\lambda|}\varepsilon t|q(t)|dt\quad(\varepsilon>0)\label{eq2.26}
\end{equation}
belong to $L^2(\mathbb{R}_+)$ it is sufficient that there exist such $\varepsilon>0$ that $e^{(1+\sqrt3)|\lambda|}\varepsilon t|q(t)|\in L^2(\mathbb{R}_+)$. Since $t\leq e^t$ ($t\geq0$), then
$$e^{(1+\sqrt3)|\lambda|}\varepsilon t|q(t)|\leq e^{[(1+\sqrt3)|\lambda|+\varepsilon]t}|q(t)|,$$
therefore, due to \eqref{eq2.3}, it is necessary to prove that one can choose $\varepsilon>0$ and $\lambda\in\mathbb{D}_{a/3}$ such that $(1+\sqrt3)|\lambda|+\varepsilon<a$ and thus
$$\varepsilon<a-(1+\sqrt3)|\lambda|+\varepsilon<a-\frac{1+\sqrt3}3a=\frac{2-\sqrt3}3a.$$
So, for all $\varepsilon$ such that ${\displaystyle0<\varepsilon<\frac{2-\sqrt3}3a}$, the function $\varepsilon c(x)$ \eqref{eq2.26} belongs to $L^2(\mathbb{R}_+)$. $\blacksquare$

{\bf Jost solutions} of system \eqref{eq2.5} are
\begin{equation}
y_k(\lambda,x)\stackrel{\rm def}{=}(a_k(\lambda)e^{i\lambda^3x},\psi_k(\lambda,x))\quad(1\leq k\leq3),\label{eq2.27}
\end{equation}
where $\lambda\mathbb{D}_{a/3}$, $\psi_k(\lambda,x)$ are given by the equations \eqref{eq2.22}, and $\{a_k(\lambda)\}$ are functions of $\lambda$. If $a_2(\lambda)=a_1(\lambda\zeta_2)$, $a_3(\lambda)=a_2(\lambda\zeta_2)$, $a_1(\lambda)=a_3(\lambda\zeta_2)$, then in view of \eqref{eq2.24} for $\{y_k(\lambda,x)\}$ the following equalities are true:
\begin{equation}
y_1(\lambda\zeta_2,x)=y_2(\lambda,x),\quad y_2(\lambda\zeta_2,x)=y_3(\lambda,x),\quad y_3(\lambda\zeta_2,x)=y_1(\lambda,x).\label{eq2.28}
\end{equation}
\vspace{5mm}

{\bf 2.2.} Find out for which $\lambda\in\mathbb{D}_{a/3}$ the functions $\{\psi_k(\lambda,x)\}_1^3$ \eqref{eq2.22} are linearly independent. For $\lambda=0$, the functions $\{\psi_k(\lambda,x)\}$ coincide, $\psi_1(\lambda,0)=\psi_2(\lambda,0)=\psi_3(\lambda,0)$. Consider the set
\begin{equation}
E_\alpha\stackrel{\rm def}{=}\{\lambda\in\mathbb{D}_{a/3}:\lambda^2(1+\alpha im_{s_2}(\lambda))=0\}.\label{eq2.29}
\end{equation}
Equation $\lambda^2(1+\alpha im_d(\lambda))=0$, due to \eqref{eq2.13}, implies that
\begin{equation}
-3i\lambda^2=\alpha\sum\zeta_km_k(\lambda)\quad(\lambda\in\mathbb{D}_{a/3})\label{eq2.30}
\end{equation}
where $\{m_k(\lambda)\}$ are given by \eqref{eq2.14}. The set $E_\alpha$ is finite, in view of analyticity of $m_d(\lambda)$ in $\mathbb{D}_{a/3}$, and $\lambda\in E_\alpha\Leftrightarrow\lambda\zeta_2\in E_\alpha$ (see \eqref{eq2.21}).

\begin{lemma}\label{l2.3}
For all $\lambda\in\mathbb{D}_{a/3}\backslash E_\alpha$ ($E_\alpha$ is given by \eqref{eq2.29}), functions $\{\psi_k(\lambda,x)\}$ \eqref{eq2.22} are linearly independent.
\end{lemma}

P r o o f. Lemma's statement follows from asymptotics \eqref{eq2.7} ($b_k(\lambda)=b(\lambda)$ \eqref{eq2.23}). Give a direct proof. Assuming the contrary, suppose that for some $\lambda\in\mathbb{D}_{a/3}$ there are such numbers $\mu_k$ ($\not=0$) from $\mathbb{C}$ ($1\leq k\leq3$) that ${\displaystyle\mu_k\psi_k(\lambda,x)=0}$ ($\forall x\in\mathbb{R}_+$), then
\begin{equation}
\sum\mu_k\psi_k(\lambda,x)=0,\quad\sum\mu_k\psi'_k(\lambda,x)=0;\quad\sum\mu_k\psi''_k(\lambda,x)=0.\label{eq2.31}
\end{equation}
Calculate the determinant (Wronskian) of this system,
$$\Delta(\lambda,x)\stackrel{\rm def}{=}\left[
\begin{array}{ccc}
\psi_1(\lambda,x)&\psi_2(\lambda,x)&\psi_3(\lambda,x)\\
\psi'_1(\lambda,x)&\psi'_2(\lambda,x)&\psi'_3(\lambda,x)\\
\psi''_1(\lambda,x)&\psi''_2(\lambda,x)&\psi''_3(\lambda,x)
\end{array}\right].$$
Using \eqref{eq2.22} and \eqref{eq1.3}, after elementary transformations, one obtains
$$\Delta(\lambda,x)=-3\sqrt3\lambda^3b^3(\lambda)-3\sqrt3ib^2(\lambda)\int\limits_x^\infty dt\int\limits_0^\infty d\tau q(t)\overline{q(\tau)}\{s_2(i\lambda(x-t))s_0(i\lambda(\tau-x)$$
$$+s_1(i\lambda(x-t))s_1(i\lambda(\tau-x))+s_0(i\lambda(x-t))s_2(i\lambda(\tau-x))\}$$
and according to (vii) \eqref{eq1.5},
$$\Delta(\lambda,x)=-3\sqrt3\lambda^3b^3(\lambda)-3\sqrt3ib^2(\lambda)\int\limits_x^\infty dt\int\limits_0^\infty d\tau s_2(i\lambda(\tau-t)q(t)\overline{q(\tau)}.$$
If $\Delta(\lambda,x)=0$ for all $x\in\mathbb{R}_+$ (i. e., solution $\{\mu_k\}_1^3$ to system \eqref{eq2.31} is non-trivial), then $\Delta'(\lambda,x)=0$ and thus
$$q(x)\int\limits_0^\infty d\tau\overline{q(\tau)}s_2(i\lambda(\tau-x))=0,$$
therefore,
$$\Delta(\lambda,x)=-3\sqrt3\lambda^3b^3(\lambda)=0$$
and thus $\lambda\in E_\alpha$ \eqref{eq2.29}. $\blacksquare$

We will need the following theorem belonging to Titchmarsh \cite{23, 25}.

\begin{theorem}[Titchmarsh]
Let $F(x)\in L^2(\mathbb{R})$, then the following statements are equivalent:

$({\rm i})$ $F(\lambda)$ is holomorphic in $\mathbb{C}_+$ and is of Hardy class $H_+^2$;

$({\rm ii})$ the first Sokhotskii formula holds,
$$\Re F(\lambda)=\frac1\pi/\hspace{-4.4mm}\int\limits_{\mathbb{R}}\frac{\Im F(x)}{x-\lambda}dx\quad(\lambda\in\mathbb{R});$$

$({\rm iii})$ the second Sokhotskii formula holds,
$$\Im F(\lambda)=-\frac1\pi/\hspace{-4.4mm}\int\limits_{\mathbb{R}}\frac{\Re F(x)}{x-\lambda}dx\quad(\lambda\in\mathbb{R}).$$
\end{theorem}

Integrals of (ii), (iii) are understood in the sense of the principal value and the formulas (ii), (iii) are called the {\bf dispersion relations} \cite{1, 4, 6, 25}.

The function $m_1^*(\lambda)$ in view of \eqref{eq2.14} is
$$m_1^*(\lambda)=\int\limits_{\mathbb{R}_+}e^{i\lambda s}\overline{g_q(s)}ds$$
and belongs to the Hardy space $H_+^2$ (Paley -- Wiener theorem \cite{24}). And since (see \eqref{eq2.17})
$$2\Re m_1^*(\lambda)=\widetilde{q}_1(\lambda)\widetilde{q}_1^*(\lambda)\quad(\lambda\in\mathbb{R}),$$
then, using (iii) of Titchmarsh theorem, one obtains
$$m_1^*(\lambda)=\frac12\widetilde{q}_1(\lambda)\widetilde{q}_1^*(\lambda)+\frac1{2\pi i}\int\limits_{\mathbb{R}}\hspace{-4.4mm}/\frac{|\widetilde{q}_1(t)|^2}{t-\lambda}dt\quad(\lambda\in\mathbb{R}),$$
and thus
\begin{equation}
m_1(\lambda)=\frac12\widetilde{q}_1(\lambda)\widetilde{q}_1^*(\lambda)-\frac1{2\pi i}\int\limits_{\mathbb{R}}\hspace{-4.4mm}/\frac{|\widetilde{q}_1(t)|^2}{t-\lambda}dt\quad(\lambda\in LD(1,a)).\label{eq2.32}
\end{equation}
Upon substitution $\lambda\rightarrow\lambda\zeta_2$, $\lambda\rightarrow\lambda\zeta_3$ in \eqref{eq2.32}, in view of \eqref{eq2.17}, one has
\begin{equation}
\begin{array}{ccc}
{\displaystyle m_2(\lambda)=\frac12\widetilde{q}_2(\lambda)\widetilde{q}_3^*(\lambda)-\frac1{2\pi i}\int\limits_{\mathbb{R}}\hspace{-4.4mm}/\frac{|\widetilde{q_1}(t)|^2}{t-\lambda\zeta_2}dt\quad(\lambda\in LD(3,a));}\\
{\displaystyle m_3(\lambda)=\frac12\widetilde{q}_3(\lambda)\widetilde{q}_2^*(\lambda)-\frac1{2\pi i}\int\limits_{\mathbb{R}}\hspace{-4.4mm}/\frac{|\widetilde{q_1}(t)|^2}{t-\lambda\zeta_3}dt\quad(\lambda\in LD(2,a)).}
\end{array}\label{eq2.33}
\end{equation}
Substituting expressions $m_k(\lambda)$ \eqref{eq2.32}, \eqref{eq2.33} into equation \eqref{eq2.30}, one obtains
$$-3i\lambda^2=\frac\alpha2\{\widetilde{q}_1(\lambda)\widetilde{q}_1^*(\lambda)+\zeta_2\widetilde{q}_2(\lambda)\widetilde{q}_3^*(\lambda)+\zeta_3\widetilde{q}_3(\lambda)\widetilde{q}_2^*(\lambda)\}$$
$$-\frac\alpha{2\pi i}\int\limits_{\mathbb{R}}\hspace{-4.3mm}/dt|\widetilde{q}_1(t)||^2\left\{\frac1{t-\lambda}+\frac{\zeta_2}{t-\lambda\zeta_2}+\frac{\zeta_3}{t-\lambda\zeta_3}\right\}\quad(\lambda\in\mathbb{D}_{a/3}\cap\mathbb{R})$$
or
$$0=\frac\alpha2\{\widetilde{q}_1(\lambda)\widetilde{q}_1^*(\lambda)+\zeta_2\widetilde{q}_2(\lambda)\widetilde{q}_3^*(\lambda)+\zeta_3\widetilde{q}_3(\lambda)\widetilde{q}_2^*(\lambda)\}+3i\lambda^2\left\{1+\frac\alpha{2\pi}\int
\hspace{-4.4mm}/\frac{|\widetilde{q}_1(t)|^2}{t^3-\lambda^3}dt\right\}.$$
The first summand in this sum is real, and the second is purely imaginary ($\forall\lambda\in\mathbb{D}_{a/3}\cap\mathbb{R}$), therefore,
$$\left\{
\begin{array}{ccc}
\widetilde{q}_1(\lambda)\widetilde{q}_1^*(\lambda)+\zeta_2\widetilde{q}_2(\lambda)\widetilde{q}_3^*(\lambda)+\zeta_3\widetilde{q}_3(\lambda)\widetilde{q}_2^*(\lambda)=0;\\
{\displaystyle\lambda^2\left\{\frac{2\pi}\alpha+\int\limits_{\mathbb{R}}\hspace{-4.4mm}/\frac{|\widetilde{q}_1(t)|^2dt}{t^3-\lambda^3}\right\}=0}
\end{array}\right.\quad(\lambda\in\mathbb{D}_{a/3}\cap\mathbb{R}).$$
The second equality, apart from $\lambda=0$, gives
\begin{equation}
\int\limits_{\mathbb{R}}\hspace{-4.4mm}/\frac{|\widetilde{q}_1(t)|^2}{t^3-\lambda^3}dt+\frac{2\pi}\alpha=0,\label{eq2.34}
\end{equation}
which (as was mentioned before) has only finite number of solutions. Since the function
$$\int\frac{|\widetilde{q}_1(t)|^2dt}{t^3-z}\quad(z\in\mathbb{C}\backslash\mathbb{R})$$
is Nevanlinna one, equation \eqref{eq2.34} does not have complex roots for $\lambda^3\not\in\mathbb{R}$.

\begin{lemma}\label{l2.4}
The set $E_\alpha$ \eqref{eq2.29} is given by
\begin{equation}
\begin{array}{ccc}
E_\alpha=\{0\}\cup\{\zeta_2^l\mu_k:\mu_k\in\mathbb{R};\widetilde{q}_1(\mu_k)\widetilde{q}_1^*(\mu_k)+\zeta_2\widetilde{q}_2(\mu_k)\widetilde{q}_3^*(\mu_k)\\
+\zeta_3\widetilde{q}_3(\mu_k)\widetilde{q}_2^*(\mu_k)=0;1\leq k\leq n<\infty,l=0,1,2\}
\end{array}\label{eq2.35}
\end{equation}
where $\mu_k$ are roots of equation \eqref{eq2.34}.
\end{lemma}

So, the set $E_\alpha$ \eqref{eq2.29} consists of a finite number ($3n$, $n\in\mathbb{N}$) points situated on the system of straight lines $L=\bigcup\limits_kL_{\zeta_k}$ \eqref{eq1.27} and this set is invariant under rotations $\lambda\rightarrow\lambda\zeta_2$, $\lambda\rightarrow\lambda\zeta_3$.

Introduce the notations that are needed further
$$\widetilde{q}_{s_0}(\lambda)\stackrel{\rm def}{=}\langle s_0(-i\lambda x),\overline{q}(x)\rangle=\frac13(\widetilde{q}_1(\lambda)+\widetilde{q}_2(\lambda)+\widetilde{q}_3(\lambda));$$
\begin{equation}
\widetilde{q}_{s_1}(\lambda)\stackrel{\rm def}{=}\langle s_1(-i\lambda x),\overline{q}(x)\rangle=\frac13(\zeta_1\widetilde{q}_1(\lambda)+\zeta_3\widetilde{q}_2(\lambda)+\zeta_2\widetilde{q}_3(\lambda));\label{eq2.37}
\end{equation}
$$\widetilde{q}_{s_2}(\lambda)\stackrel{\rm def}{=}\langle s_2(-i\lambda x),\overline{q(x)}\rangle=\frac13(\zeta_1\widetilde{q}_1(\lambda)+\zeta_2\widetilde{q}_2(\lambda)+\zeta_3\widetilde{q}_3(\lambda)).$$

\section{Scattering problem}\label{s3}

{\bf 3.1.} Linear independence of the functions $\{\psi_k(\lambda,x)\}_1^3$ for $\lambda\in\mathbb{D}_{a/3}\backslash E_\alpha$ (Lemma \ref{l2.3}) implies that Jost solutions $\{y_k(\lambda,x)\}$ \eqref{eq2.27} to system \eqref{eq2.5} also are linearly independent for these $\lambda$, therefore an arbitrary solution $y(\lambda,x)$ to the equations system \eqref{eq2.5} is their linear combination
\begin{equation}
y(\lambda,x)=\sum\limits_kB_ky_k(\lambda,x)=\left(\sum a_kB_ke^{i\lambda^3x},\sum\limits_kB_k\psi_k(\lambda,x)\right)\label{eq3.1}
\end{equation}
where $\{a_k\}$, $\{B_k\}$ do not depend on $x$. Excluding the trivial case ($B_k=0$, $\forall k$), suppose that $B_1$ (for example) is nonzero, $B_1\not=0$, then
\begin{equation}
\begin{array}{ccc}
y(\lambda,x)B_1^{-1}=y_1(\lambda,x)+S_2y_2(\lambda,x)+S_3y_3(\lambda,x)\\
=(a(\lambda)e^{i\lambda^3x},\psi_1(\lambda,x)+S_2(\lambda)\psi_2(\lambda,x)+S_3(\lambda)\psi_3(\lambda,x)),
\end{array}\label{eq3.2}
\end{equation}
here
$$a=B_1^{-1}\left(\sum B_ka_k\right),\quad S_2=B_1^{-1}B_2,\quad S_3=B_1^{-1}B_3.$$

\begin{remark}\label{r3.1}
Using \eqref{eq2.7}, one obtains that
\begin{equation}
y(\lambda,x)B_1^{-1}\rightarrow b(\lambda)(e^{i\lambda\zeta_1x}+S_2(\lambda)e^{i\lambda\zeta_2x}+S_3(\lambda)e^{i\lambda\zeta_3x})\quad(x\rightarrow\infty)\label{eq3.3}
\end{equation}
where $b(\lambda)$ is given by \eqref{eq2.23} and thus $S_2(\lambda)$ and $S_3(\lambda)$ should be considered as the {\bf scattering coefficients} of the incident wave $e^{i\lambda\zeta_1x}$. The function $a(\lambda)$ is said to be the {\bf matching coefficient} of the incident wave $e^{i\lambda\zeta_1x}$. Equation \eqref{eq3.3}, after the substitutions $\lambda\rightarrow\lambda\zeta_2$, $\lambda\rightarrow\lambda\zeta_3$, easily yields scattering coefficients of incident waves $e^{i\lambda\zeta_2x}$ and $e^{i\lambda\zeta_3x}$.
\end{remark}

In order to find the coefficients $S_2(\lambda)$ and $S_3(\lambda)$, use boundary conditions \eqref{eq2.4} for the function $y(\lambda,x)B^{-1}$ \eqref{eq3.2},
\begin{equation}
\left\{
\begin{array}{lll}
S_2(\lambda)\psi_2(\lambda,0)+S_3(\lambda)\psi_3(\lambda,0)=-\psi_1(\lambda,0);\\
S_2(\lambda)\psi'_2(\lambda,0)+S_3(\lambda)\psi'_3(\lambda,0)=-\psi'_1(\lambda,0)+a(\lambda).
\end{array}\right.\label{eq3.4}
\end{equation}
Determinant $\Delta(\lambda)$ of this system equals to the Wronskian $\Delta(\lambda)=W_{2,3}(\lambda)$ where
\begin{equation}
W_{k,s}(\lambda)\stackrel{\rm def}{=}\psi_k(\lambda,0)\psi'_s(\lambda,0)-\psi_s(\lambda,0)\psi'_k(\lambda,0)\quad(1\leq k,s\leq3).\label{eq3.5}
\end{equation}
Evidently, $W_{k,s}(\lambda)\not=0$ ($k\not=s$) for all $\lambda\in\mathbb{D}_{a/3}\backslash E_{\alpha}$; thus, $W_{k,s}(\lambda)=0$ implies (see \eqref{eq3.5}) linear dependance of $\psi_k(\lambda,0)$ and $\psi_s(\lambda,0)$ which contradicts Lemma \ref{l2.3}. From equation \eqref{eq3.4} one finds that for all $\lambda\in\mathbb{D}_{a/3}\backslash E_\alpha$
\begin{equation}
S_2(\lambda)=\frac{W_{3,1}(\lambda)-a(\lambda)\psi_3(\lambda,0)}{W_{2,3}(\lambda)};\quad S_3(\lambda)=\frac{W_{1,2}(\lambda)+a(\lambda)\psi_2(\lambda,0)}{W_{2,3}(\lambda)}.\label{eq3.6}
\end{equation}

\begin{lemma}\label{l3.1}
The function $W_{2,3}(\lambda)$ \eqref{eq3.5} has representation
\begin{equation}
W_{2,3}(\lambda)=\zeta_1\sqrt3\lambda b(\lambda)\psi_1^*(\lambda,0)\quad(\forall\lambda\in\mathbb{D}_{a/3}\backslash E_\alpha)\label{eq3.7}
\end{equation}
where $b(\lambda)$ is given by \eqref{eq2.23}.
\end{lemma}

P r o o f. Equations \eqref{eq2.22} and \eqref{eq2.37} imply
$$W_{2,3}(\lambda)=\left(b(\lambda)+\frac{\alpha i}{\lambda^2}\widetilde{q}_3^*(\lambda)\widetilde{q}_{s_2}(\lambda)\right)\left(i\lambda\zeta_3b(\lambda)-\frac\alpha\lambda\widetilde{q}_2^*(\lambda)\widetilde{q}_{s_1}(\lambda)\right)$$
$$-\left(b(\lambda)+\frac{\alpha i}{\lambda^2}\widetilde{q}_2^*(\lambda)\widetilde{q}_{s_2}(\lambda)\right)\left(i\lambda\zeta_2b(\lambda)-\frac\alpha\lambda\widetilde{q}_3^*(\lambda)\widetilde{q}_{s_1}(\lambda\right)=i\lambda(\zeta_3-\zeta_2)b^2(\lambda)$$
$$-\frac{\alpha b(\lambda)}\lambda\{\widetilde{q}_2^*(\lambda)\langle[s_1(-i\lambda x)-\zeta_2s_2(-i\lambda x)],\overline{q}(x)\rangle-\widetilde{q}_3^*(\lambda)\langle[s_1(-i\lambda x)-\zeta_3s_2(-i\lambda x)],\overline{q}(x)\rangle\}.$$
Since, in view of \eqref{eq1.3},
$$s_1(z)-\zeta_2s_2(z)=\frac13(\zeta_1-\zeta_2)\left(e^{z\zeta_1}-e^{z\zeta_3}\right);\quad s_1(z)-\zeta_3d(z)=\frac13(\zeta_1-\zeta_3)\left(e^{z\zeta_1}-e^{z\zeta_2}\right),$$
then
$$W_{2,3}(\lambda)=\sqrt3\lambda b^2(\lambda)-\frac{\alpha b(\lambda)}\lambda\{\zeta_3\widetilde{q}_2^*(\lambda)[\widetilde{q}_1(\lambda)-\widetilde{q}_3(\lambda)]+\zeta_2\widetilde{q}_3^*(\lambda)[\widetilde{q}_1(\lambda)-\widetilde{q}_2(\lambda)]\}$$
$$=\frac{\sqrt3b(\lambda)}\lambda\left\{\lambda^2-\frac{\alpha i}3\sum\zeta_km_k(\lambda)-\frac{\alpha i}3[\widetilde{q}_1(\lambda)(\zeta_3\widetilde{q}_2^*(\lambda)+\zeta_2\widetilde{q}_3^*(\lambda))-\zeta_3\widetilde{q}_3(\lambda)\widetilde{q}_2^*(\lambda)\right.$$
$$-\zeta_2\widetilde{q}_2(\lambda)\widetilde{q}_3^*(\lambda)\}$$
and due to \eqref{eq2.17} one obtains
$$W_{2,3}(\lambda)=\frac{\sqrt3b(\lambda)}\lambda\left\{\lambda^2+\frac{\alpha i}3\left[\left(\sum\zeta_km_k(\lambda)\right)^*-\widetilde{q}_1(\lambda)\left(\sum\zeta_k\widetilde{q}_k(\lambda)\right)^*\right]\right\}$$
$$=\sqrt3b(\lambda)\left[b^*(\lambda)-\frac{\alpha i}{\lambda^2}\widetilde{q}_1(\lambda)\widetilde{q}_d^*(\lambda)\right].\blacksquare$$

\begin{corollary}\label{c3.1}
Equations \eqref{eq2.24} and $\psi'_1(\lambda\zeta_2,x)=\psi'_2(\lambda,x)$, $\psi'_2(\lambda\zeta_2,x)=\psi'_3(\lambda,x)$, $\psi'_3(\lambda\zeta_2,x)=\psi'_1(\lambda,x)$ imply that
$$W_{2,3}(\lambda\zeta_2)=W_{3,1}(\lambda);\quad W_{3,1}(\lambda\zeta_2)=W_{1,2}(\lambda),$$
hence, in view of \eqref{eq3.7},
\begin{equation}
W_{3,1}(\lambda)=\zeta_2\sqrt3\lambda b(\lambda)\psi_3^*(\lambda,0);\quad W_{2,1}(\lambda)=\zeta_3\sqrt3\lambda b(\lambda)\psi_2^*(\lambda)\quad(\forall\lambda\in\mathbb{D}_{a/3}\backslash E_\alpha).\label{eq3.8}
\end{equation}
The following relations hold:
\begin{equation}
\begin{array}{lll}
{\displaystyle\sum\limits_k\zeta_k\psi_k^*(\lambda,0)=-\frac{3\alpha i}{\lambda^2}q_{s_2}(\lambda)\widetilde{q}_{s_2}^*(\lambda);}\\
{\displaystyle W_{2,3}(\lambda)+W_{3,1}(\lambda)+W_{1,2}(\lambda)=-\frac{3\sqrt3}\lambda b(\lambda)\widetilde{q}_{s_2}(\lambda)\widetilde{q}_{s_2}^*(\lambda)}.
\end{array}\label{eq3.9}
\end{equation}
\end{corollary}
\vspace{5mm}

{\bf 3.2.} Next, study properties of scattering coefficients $S_2(\lambda)$ and $S_3(\lambda)$ \eqref{eq3.6}. Equations \eqref{eq3.6}, taking into account \eqref{eq3.7}, \eqref{eq3.8}, give system of equations
\begin{equation}
\left\{
\begin{array}{lll}
\zeta_1\sqrt3\lambda b(\lambda)\psi_1^*(\lambda,0)S_2(\lambda)=\zeta_2\sqrt3\lambda b(\lambda)\psi_3^*(\lambda,0)-a(\lambda)\psi_3(\lambda,0);\\
\zeta_1\sqrt3\lambda b(\lambda)\psi_1^*(\lambda,0)S_3(\lambda)=\zeta_3\sqrt3\lambda b(\lambda)\psi_2^*(\lambda,0)+a(\lambda)\psi_2(\lambda,0)
\end{array}\right.\label{eq3.10}
\end{equation}
where $\lambda\in\mathbb{D}_{a/3}\backslash E_\alpha$. Making the change of variables $\lambda\rightarrow\lambda\zeta_2$ in the second equation of this system, one obtains
$$\zeta_2\sqrt3\lambda b(\lambda)\psi_3^*(\lambda,0)S_3(\lambda\zeta_2)=\zeta_1\sqrt3\lambda b(\lambda)\psi_1^*(\lambda,0)+a(\lambda\zeta_2)\psi_3(\lambda,0).$$
Upon multiplying this equation by $a(\lambda)$ and the first equation in \eqref{eq3.10}, correspondingly, by $a(\lambda\zeta_2)$, and adding, one obtains the relation
$$\zeta_2\psi_3^*(\lambda,0)[S_3(\lambda\zeta_2)a(\lambda)-a(\lambda\zeta_2)]=\zeta_1\psi_1^*(\lambda,0)[a(\lambda)-S_2(\lambda)a(\lambda\zeta_2)],$$
and thus
$$\psi_1^*(\lambda,0)=-\zeta_2T(\lambda)\psi_3^*(\lambda,0)$$
where
\begin{equation}
T(\lambda)\stackrel{\rm def}{=}\frac{S_3(\lambda\zeta_2)a(\lambda)-a(\lambda\zeta_2)}{S_2(\lambda)a(\lambda\zeta_2)-a(\lambda)}.\label{eq3.11}
\end{equation}

\begin{lemma}\label{l3.2}
For all $\lambda\in\mathbb{D}_{a/3}\backslash E_\alpha$, the following equations are true:
\begin{equation}
\begin{array}{ccc}
\psi_1^*(\lambda,0)=-\zeta_2T(\lambda)\psi_3^*(\lambda,0);\quad\psi_3^*(\lambda,0)=-\zeta_2T(\lambda\zeta_2)\psi_2^*(\lambda,0);\\
\psi_2^*(\lambda,0)=-\zeta_2T(\lambda\zeta_3)\psi_1^*(\lambda,0)
\end{array}\label{eq3.12}
\end{equation}
where $T(\lambda)$ is given by \eqref{eq3.11} and
\begin{equation}
T(\lambda)T(\lambda\zeta_2)T(\lambda\zeta_3)=-1.\label{eq3.13}
\end{equation}
\end{lemma}

Thus, knowing one of the functions $\psi_k^*(\lambda,0)$ and using \eqref{eq3.12}, one can find the other functions $\psi_p^*(\lambda,0)$ ($p\not=k$).

Upon substitution of $\psi_1^*(\lambda,0)$ from the first equation of \eqref{eq3.12} into the first equation of system \eqref{eq3.10}, one obtains
$$\zeta_2\sqrt3b(\lambda)\psi_3^*(\lambda,0)(S_2(\lambda)T(\lambda)+1)=a(\lambda)\psi_3(\lambda,0)$$
or
$$\psi_3^*(\lambda,0)=u(\lambda)\psi_3(\lambda,0)$$
where
\begin{equation}
u(\lambda)=\frac{\zeta_3a(\lambda)}{\sqrt3b(\lambda)[S_2(\lambda)T(\lambda)+1]}.\label{eq3.14}
\end{equation}

\begin{lemma}\label{l3.3}
For all $\lambda\in\mathbb{D}_{a/3}\backslash E_\alpha$,
\begin{equation}
\begin{array}{ccc}
\psi_3^*(\lambda,0)=u(\lambda)\psi_3(\lambda,0);\quad\psi_2^*(\lambda,0)=u(\lambda\zeta_2)\psi_1(\lambda,0);\\
\psi_1^*(\lambda,0)=u(\lambda\zeta_3)\psi_2(\lambda,0)
\end{array}\label{eq3.15}
\end{equation}
where $u(\lambda)$ is given by \eqref{eq3.14} and
\begin{equation}
u^*(\lambda)u(\lambda)=1;\quad u^*(\lambda\zeta_3)u(\lambda\zeta_2)=1.\label{eq3.16}
\end{equation}
\end{lemma}

The following statement contains a description of properties of the functions $\{\psi_k(\lambda,0)\}$.

\begin{lemma}\label{l3.4}
For the functions $\{\psi_k(\lambda,0)\}$ for all $\lambda\in\mathbb{D}_{a/3}$ the following equalities hold:

\begin{equation}
\begin{array}{lll}
({\rm i})&{\displaystyle\psi_3^*(\lambda,0)-\psi_3(\lambda,0)=\frac{\alpha i}{3\lambda^2}(\widetilde{q}_1(\lambda)-\widetilde{q}_2(\lambda))(\widetilde{q}_1^*(\lambda)-\widetilde{q}_2^*(\lambda));}\\
({\rm ii})&{\displaystyle\psi_2^*(\lambda,0)-\psi_2(\lambda,0)=\frac{\alpha i}{3\lambda^2}(\widetilde{q}_1(\lambda)-\widetilde{q}_3(\lambda))(\widetilde{q}_1^*(\lambda)-\widetilde{q}_3^*(\lambda));}\\
({\rm iii})&{\displaystyle\psi_1^*(\lambda)-\psi_1(\lambda,0)=\frac{\alpha i}{3\lambda^2}\{\zeta_2(\widetilde{q}_1(\lambda)-\widetilde{q}_2(\lambda))(\widetilde{q}_1^*(\lambda)-\widetilde{q}_3^*(\lambda))}\\
&-\zeta_3(\widetilde{q}_1(\lambda)-\widetilde{q}_3(\lambda))(\widetilde{q}_1^*(\lambda)-\widetilde{q}_2^*(\lambda))\}
\end{array}\label{eq3.17}
\end{equation}
where $\{\widetilde{q}_k(\lambda)\}$ is given by \eqref{eq2.12}.
\end{lemma}

P r o o f. First notice that \eqref{eq2.20} implies
\begin{equation}
\begin{array}{ccc}
{\displaystyle b(\lambda)-b^*(\lambda)=\alpha i(m_{s_2}(\lambda)-m_{s_2}^*(\lambda))=\frac{\alpha i}{3\lambda^2}\{\widetilde{q}_1(\lambda)\widetilde{q}_1^*(\lambda)+\zeta_2\widetilde{q}_2(\lambda)\widetilde{q}_3^*(\lambda)}\\
+\zeta_3\widetilde{q}_3(\lambda)\widetilde{q}_2^*(\lambda)\}.
\end{array}\label{eq3.18}
\end{equation}
Using \eqref{eq3.18}, one finds that
$$\psi_3^*(\lambda,0)-\psi_3(\lambda,0)=b^*(\lambda)-b(\lambda)-\frac{\alpha i}{\lambda^2}[\widetilde{q}_2(\lambda)\widetilde{q}_{s_2}^*(\lambda)+\widetilde{q}_{s_2}(\lambda)\widetilde{q}_2^*(\lambda)]$$
$$=\frac{\alpha i}{3\lambda^2}\{\widetilde{q}_1(\lambda)\widetilde{q}_1^*(\lambda)+\zeta_2\widetilde{q}_2(\lambda)\widetilde{q}_3^*(\lambda)+\zeta_3\widetilde{q}_3(\lambda)\widetilde{q}_2^*(\lambda)-\widetilde{q}_2(\lambda)[\widetilde{q}_1^*
(\lambda)+
\zeta_3\widetilde{q}_2^*(\lambda)+\zeta_2\widetilde{q}_3^*(\lambda)]$$
$$-[\widetilde{q}_1(\lambda)+\zeta_2\widetilde{q}_2(\lambda)+\zeta_3\widetilde{q}_3(\lambda)]\}=\frac{\alpha i}{3\lambda^2}\{q_1(\lambda)\widetilde{q}_1^*(\lambda)-\widetilde{q}_2(\lambda)\widetilde{q}_1^*(\lambda)-\widetilde{q}_1(\lambda)\widetilde{q}_2^*(\lambda)$$
$$-(\zeta_2+\zeta_3)\widetilde{q}_2(\lambda)\widetilde{q}_2^*(\lambda)\},$$
hence it follows (i) \eqref{eq3.17}, due to $\zeta_2+\zeta_3=-\zeta_1$ ($=-1$). Equation (ii) \eqref{eq3.17} is proved analogously. In order to prove equation (iii) \eqref{eq3.17}, again use \eqref{eq3.18}, then
$$\psi_1^*(\lambda,0)-\psi_1(\lambda,0)=\frac{\alpha i}{3\lambda^2}\{q_1(\lambda)\widetilde{q}_1^*(\lambda)+\zeta_2\widetilde{q}_2(\lambda)\widetilde{q}_3^*(\lambda)+\zeta_3\widetilde{q}_3(\lambda)\widetilde{q}_2^*(\lambda)$$
$$-\widetilde{q}_1(\lambda)[\widetilde{q}_1^*(\lambda)+\zeta_3\widetilde{q}_2^*(\lambda)+\zeta_2\widetilde{q}_3^*(\lambda)]-[\widetilde{q}_1(\lambda)+\zeta_2\widetilde{q}_2(\lambda)+\zeta_3\widetilde{q}_3(\lambda)]\widetilde{q}_1^*(
\lambda)\}$$
$$=\frac{\alpha i}{3\lambda^2}\{-\widetilde{q}_1(\lambda)\widetilde{q}_1^*(\lambda)+\zeta_2\widetilde{q}_2(\lambda)(\widetilde{q}_3^*(\lambda)-\widetilde{q}_1^*(\lambda))+\zeta_3\widetilde{q}_3(\lambda)(\widetilde{q}_2^*(\lambda)-\widetilde{q}_1^*
(\lambda))$$
$$-\zeta_3\widetilde{q}_1(\lambda)\widetilde{q}_2^*(\lambda)-\zeta_2\widetilde{q}_1(\lambda)\widetilde{q}_3^*(\lambda)\},$$
this gives (iii) \eqref{eq3.17} since $\zeta_2+\zeta_3=-1$. $\blacksquare$

\begin{corollary}\label{c3.2}
For all $\lambda\in\mathbb{D}_{a/3}$, the following equality holds,
\begin{equation}
\zeta_3\psi_2^*(\lambda,0)+\zeta_2\psi_3^*(\lambda,0)+\zeta_1\psi_1(\lambda,0)=-\frac{\alpha i}{3\lambda^2}[\widetilde{q}_2(\lambda)-\widetilde{q}_3(\lambda)][\widetilde{q}_2^*(\lambda)-\widetilde{q}_3^*(\lambda)].\label{eq3.19}
\end{equation}
\end{corollary}

P r o o f. Substitute $\lambda\rightarrow\lambda\zeta_2$ into equation (i) \eqref{eq3.17}, then
$$\psi_2^*(\lambda,0)-\psi_1(\lambda,0)=\frac{\alpha i\zeta_2}{3\lambda^2}(\widetilde{q}_2(\lambda)-\widetilde{q}_3(\lambda))(\widetilde{q}_3^*(\lambda)-\widetilde{q}_1^*(\lambda)).$$
Similarly, after the substitution $\lambda\rightarrow\lambda\zeta_3$ in (ii) \eqref{eq3.17}, one obtains
$$\psi_3^*(\lambda,0)-\psi_1(\lambda,0)=\frac{\alpha i\zeta_3}{3\lambda^2}(\widetilde{q}_3(\lambda)-\widetilde{q}_2(\lambda))(\widetilde{q}_2^*(\lambda)-\widetilde{q}_1^*(\lambda)).$$
Upon multiplying the first of the obtained equations by $\zeta_3$ and the second, correspondingly, by $\zeta_2$ and adding, one has relation \eqref{eq3.19}. $\blacksquare$
\vspace{5mm}

{\bf 3.3.} Taking \eqref{eq2.17} into account, rewrite $\psi_1(\lambda,0)$ as
\begin{equation}
\psi_1(\lambda,0)=1-\frac{\alpha i}{3\lambda^2}\{-m_1^*(\lambda)+\zeta_2(m_2(\lambda)-\widetilde{q}_2(\lambda)\widetilde{q}_1^*(\lambda))+\zeta_3(m_3(\lambda)-q_3(\lambda)\widetilde{q}_1^*(\lambda))\}.\label{eq3.20}
\end{equation}
Hence (see Subseq. 1.2) follows holomorphy of $\psi_1(\lambda,0)$ in the sector $S_2$ (see \eqref{eq1.11}).

\begin{picture}(200,200)
\put(0,100){\vector(1,0){200}}
\put(150,0){\vector(-1,2){100}}
\put(150,200){\vector(-1,-2){100}}
\put(30,190){$L_{\zeta_2}$}
\put(69,190){$\psi_1(\lambda,0)$}
\put(140,150){$\psi_3^*(\lambda,0)$}
\put(20,150){$\psi_2^*(\lambda,0)$}
\put(20,50){$\psi_3(\lambda,0)$}
\put(80,108){$S_3$}
\put(80,88){$S_4$}
\put(92,115){$S_2$}
\put(90,70){$S_5$}
\put(115,87){$S_6$}
\put(105,103){$S_1$}
\put(190,105){$L_{\zeta_1}$}
\put(130,60){$\psi_2(\lambda,0)$}
\put(80,40){$\psi_1^*(\lambda,0)$}
\put(60,0){$L_{\zeta_3}$}
\end{picture}

\hspace{20mm} Fig. 4

Since $\psi_2(\lambda,0)=\psi_1(\lambda\zeta_2,0)$ (see \eqref{eq2.24}), $\psi_2(\lambda,0)$ is analytical in $S_6$. Similarly, $\psi_1(\lambda\zeta_3,0)=\psi_3(\lambda,0)$ implies holomorphy of the function $\psi_3(\lambda,0)$ in the sector $S_4$. This provides analyticity of $\psi_1^*(\lambda,0)$, $\psi_2^*(\lambda,0)$, $\psi_3^*(\lambda,0)$ in the sectors $S_5$, $S_3$, $S_1$ correspondingly (see Fig. 4).

Equations \eqref{eq3.12}, \eqref{eq3.15} for the functions $\{\psi_k(\lambda,0)\}$ and $\{\psi_k^*(\lambda,0)\}$ on the contour $L=\bigcup\limits_k L_{\zeta_k}$ \eqref{eq1.23} formed by the straight lines $L_{\zeta_k}$ \eqref{eq1.9} imply a boundary value problem \cite{26, 27}. So, on the straight line $L_{\zeta_2}$ boundary conditions are
\begin{equation}
\begin{array}{ccc}
\psi_2^*(\lambda,0)=U(\lambda\zeta_2)\psi_1(\lambda,0)&(\lambda\in l_{\zeta_2})\\
\psi_1^*(\lambda,0)=U(\lambda\zeta_3)\psi_2(\lambda,0)&(\lambda\in\widehat{l}_{\zeta_2})
\end{array}\label{eq3.21}
\end{equation}
on the straight line $L_{\zeta_3}$,
\begin{equation}
\begin{array}{ccc}
\psi_1^*(\lambda,0)=-\zeta_2U(\lambda)T(\lambda)\psi_3(\lambda,0)&(\lambda\in l_{\zeta_3});\\
\psi_3^*(\lambda,0)=-\zeta_2U(\lambda\zeta_2)T(\lambda\zeta_2\psi_1(\lambda,0)&(\lambda\in\widehat{l}_{\zeta_3});
\end{array}\label{eq3.22}
\end{equation}
and, finally, on $L_{\zeta_1}$,
\begin{equation}
\begin{array}{ccc}
\psi_3^*(\lambda,0)=-\zeta_3U(\lambda\zeta_3)T^{-1}(\lambda)\psi(\lambda,0)&(\lambda\in l_{\zeta_1});\\
\psi_2^*(\lambda,0)=-\zeta_3U(\lambda)T^{-1}(\lambda\zeta_2)\psi_3(\lambda,0)&(\lambda\in\widehat{l}_{\zeta_1}).
\end{array}\label{eq3.23}
\end{equation}
Coefficients of boundary value problems \eqref{eq3.21} -- \eqref{eq3.23} initially given for $\lambda\in\mathbb{D}_{a/3}$ are naturally defined on the rays $\{l_{\zeta_k}\}$ and $\{\widehat{l}_{\zeta_k}\}$ \eqref{eq1.10} in view of analyticity of the functions $\{\psi_k(\lambda,0)\}$ and $\{\psi_k^*(\lambda,0)\}$ inside of the corresponding sectors $S_p$ \eqref{eq1.11} (see $LD(k,a)$ \eqref{eq2.16}).

\begin{lemma}\label{l3.5}
Coefficients of the boundary value problems \eqref{eq3.21} -- \eqref{eq3.23}
\begin{equation}
\begin{array}{ccc}
G_2(\lambda)\stackrel{\rm def}{=}U(\lambda\zeta_2)\,(\lambda\in l_{\zeta_2});\quad G_3(\lambda)\stackrel{\rm def}{=}-\zeta_2U(\lambda)T(\lambda)\,(\lambda\in l_{\zeta_3});\\
G_1(\lambda)\stackrel{\rm def}{=}-\zeta_3U(\lambda\zeta_3)T^{-1}(\lambda)\,(\lambda\in l_{\zeta_1})
\end{array}\label{eq3.24}
\end{equation}
($U(\lambda)$ and $T(\lambda)$ are given by \eqref{eq3.11}, \eqref{eq3.14}) satisfy relations
\begin{equation}
\begin{array}{ccc}
G_2(\lambda)G_2^*(-\lambda\zeta_2)=1\,(\forall\lambda\in l_{\zeta_2});\quad G_3(\lambda)G_3^*(-\lambda\zeta_2)=1\,(\forall\lambda\in l_{\zeta_3});\\
G_1(\lambda)G_1^*(-\lambda\zeta_2)=1\,(\forall\lambda\in l_{\zeta_1}).
\end{array}\label{eq3.25}
\end{equation}
\end{lemma}

P r o o f. The second equation in \eqref{eq3.21} implies that
$$\psi_1(\lambda,0)=U^*(\lambda\zeta_3)\psi_2^*(\lambda,0)\quad(\lambda\in\widehat{l}_{\zeta_2}).$$
And since mapping $\lambda\rightarrow-\lambda$ maps the semi-axis $\widehat{l}_{\zeta_2}$ into $l_{\zeta_2}$, then from the first equation in \eqref{eq3.21} one has $U(\lambda\zeta_2)U^*(-\lambda\zeta_3)=1$ for all $\lambda\in l_{\zeta_2}$, this proves the first equation in \eqref{eq3.26}. Other equations in \eqref{eq3.25} are proved analogously. $\blacksquare$

\begin{remark}
One ought to consider equations \eqref{eq3.25} as an analogue of the unitarity condition for a scattering matrix \cite{1, 3, 4, 6, 25}. Moreover, equations \eqref{eq3.25} give rule for extending the coefficients $\{G_k(\lambda)\}$ \eqref{eq3.24} onto the semi-axes $\{\widehat{l}_{\zeta_k}\}$ \eqref{eq1.10}.

Holomorphic in the sector $S_2$ function $\psi_1(\lambda,0)$ \eqref{eq3.20} has asymptotic $\psi_1(\lambda,0)\rightarrow1$ for $\lambda\rightarrow\infty$ (for all $\lambda\in S_2$). Hence (after the substitutions $\lambda\rightarrow\lambda\zeta_2$ and $f(\lambda)\rightarrow f^*(\lambda)$) it follows that the functions $\{\psi_k(\lambda,0)\}$ and $\{\psi_k^*(\lambda,0)\}$ have the same property inside the corresponding sectors $S_p$ \eqref{eq1.11}.
\end{remark}
\vspace{5mm}

{\bf 3.4.} Consider the Jost solution \eqref{eq2.25} to system \eqref{eq2.5}
\begin{equation}
y_1(\lambda,x)=(a(\lambda)e^{i\lambda^3x},\psi_1(\lambda,x))\label{eq3.26}
\end{equation}
where $\psi_1(\lambda,x)$ is given by \eqref{eq2.22}.

\begin{remark}\label{r3.5}
The component $a(\lambda)e^{i\lambda^3x}$ in \eqref{eq3.26} for all $\lambda\in\mathbb{D}_{a/3}\cap S_2$ belongs to $L^2(\mathbb{R}_-)$ and the second component $\psi_1(\lambda,x)$ belongs to $L^2(\mathbb{R}_+)$ for all $\lambda\in\mathbb{D}_{a/3}\cap\mathbb{C}_+$ (Lemma \ref{l2.2}). So, if $q$ satisfies relation \eqref{eq2.3}, then the function $y_1(\lambda,x)$ \eqref{eq3.26} belongs to the space $\mathcal{H}$ \eqref{eq2.1} for all $\lambda\in\mathbb{D}_{a/3}\cap S_2$.
\end{remark}

Find out for which $\lambda\in\mathbb{D}_{a/3}\cap S_2$ the function $y_1(\lambda,x)$ \eqref{eq3.26} belongs to domain $\mathfrak{D}(\mathcal{L}_\alpha)$ of the operator $\mathcal{L}_\alpha$. Equation \eqref{eq2.4} for $y_1(\lambda,x)$ \eqref{eq3.26} yields system of equations
\begin{equation}
\left\{
\begin{array}{lll}
{\displaystyle b(\lambda)+\frac{\alpha i}{\lambda^2}\widetilde{q}_1^*(\lambda)\int\limits_0^\infty s_2(-i\lambda t)q(t)dt=0;}\\
{\displaystyle i\lambda b(\lambda)-\frac\alpha\lambda\widetilde{q}_1^*(\lambda)\int\limits_0^\infty s_1(-i\lambda t)q(t)dt=a(\lambda).}
\end{array}\right.\label{eq3.27}
\end{equation}

\begin{lemma}\label{l3.6}
The system of equations \eqref{eq3.27} is solvable then and only then when $\lambda\in E_\alpha$ \eqref{eq2.29} and $a(\lambda)=0$.
\end{lemma}

P r o o f. The first equation in \eqref{eq3.27} means that $\psi_1(\lambda,0)=0$ and thus $\psi_2(\lambda,0)=0$, $\psi_3(\lambda,0)=0$, in view of \eqref{eq2.24}. Using the identity
$$\zeta_1\psi_1(\lambda,0)+\zeta_3\psi_2(\lambda,0)+\zeta_2\psi_3(\lambda,0)=\frac{3i\alpha}{\lambda^2}\widetilde{q}_{s_2}(\lambda)\widetilde{q}_{s_2}^*(\lambda)$$
following from \eqref{eq3.9}, one obtains that
\begin{equation}
\widetilde{q}_d(\lambda)=\frac13(\widetilde{q}_1(\lambda)+\zeta_2\widetilde{q}_2(\lambda)+\zeta_3\widetilde{q}_3(\lambda))=0,\label{eq3.28}
\end{equation}
hence it follows (see the first equation in \eqref{eq3.27}) that $b(\lambda)=0$, i. e., $\lambda\in E_\alpha$ \eqref{eq2.35} (Lemma \ref{l2.4}). Using \eqref{eq3.28}, one finds
$$\widetilde{q}_1(\lambda)+\zeta_3\widetilde{q}_2(\lambda)+\zeta_2\widetilde{q}_3(\lambda)=(\zeta_3-\zeta_2)(\widetilde{q}_2(\lambda)-\widetilde{q}_3(\lambda))=i\sqrt3(\widetilde{q}_2(\lambda)-\widetilde{q}_3(\lambda));$$
\begin{equation}
\widetilde{q}_1(\lambda)-\widetilde{q}_2(\lambda)=\zeta_3(\widetilde{q}_2(\lambda)-\widetilde{q}_3(\lambda));\label{eq3.29}
\end{equation}
$$\widetilde{q}_1(\lambda)-\widetilde{q}_3(\lambda)=\zeta_2(\widetilde{q}_3(\lambda)-\widetilde{q}_2(\lambda));$$
and thus ($b(\lambda)=0$),
$$\psi_2(\lambda,0)=\frac{\alpha\sqrt3}{3\lambda^2}\widetilde{q}_3^*(\lambda)(\widetilde{q}_2(\lambda)-\widetilde{q}_3(\lambda))=0;$$
$$\psi_3(\lambda,0)=\frac{\alpha\sqrt3}{3\lambda^2}\widetilde{q}_2^*(\lambda)(\widetilde{q}_2(\lambda)-\widetilde{q}_3(\lambda))=0.$$
Subtracting these equations, one obtains that
$$(\widetilde{q}_2(\lambda)-\widetilde{q}_3(\lambda))(\widetilde{q}_2(\lambda)-\widetilde{q}_3(\lambda))^*=0,$$
i. e., $\widetilde{q}_2(\lambda)=\widetilde{q}_3(\lambda)$ and in view of \eqref{eq3.29} $\widetilde{q}_1(\lambda)=\widetilde{q}_2(\lambda)$. And since
$$\int\limits_0^\infty s_1(-i\lambda t)q(t)dt=\frac13(\widetilde{q}_1(\lambda)+\zeta_3\widetilde{q}_2(\lambda)+\zeta_2\widetilde{q}_3(\lambda)),$$
then $q_1(\lambda)=\widetilde{q}_2(\lambda)=\widetilde{q}_3(\lambda)$ implies that this expression vanishes, therefore the second equation in \eqref{eq3.27} yields $a(\lambda)=0$. So, necessity of the lemma's conditions is proved. Sufficiency of the statement is obvious. $\blacksquare$

\begin{remark}\label{r3.3}
Let the matching coefficient be
\begin{equation}
a(\lambda)=C(\lambda)b(\lambda)\label{eq3.30}
\end{equation}
where $b(\lambda)$ is given by \eqref{eq2.23} and $C(\lambda)$ is a function of $\lambda$ (such that \eqref{eq3.30} for coefficient $a(\lambda)$ follows from the second boundary condition \eqref{eq2.4}). Then the second condition of Lemma \ref{l3.6} follows from the first one. In this case, due to $b(\lambda\zeta_2)=b(\lambda)$, the functions $T(\lambda)$ \eqref{eq3.11} and $U(\lambda)$ \eqref{eq3.14} are
\begin{equation}
T(\lambda)=\frac{S_3(\lambda\zeta_2)C(\lambda)-C(\lambda\zeta_2)}{S_2(\lambda)C(\lambda\zeta_2)-C(\lambda)};\quad U(\lambda)=\frac1{\sqrt{3}}\zeta_3C(\lambda)\frac{1-S_2^{-1}(\lambda)}{S_3(\lambda\zeta_2)-1},\label{eq3.31}
\end{equation}
and if $C(\zeta_2\lambda)=\zeta_2^rC(\lambda)$ ($r=0$, $1$, $2$), then
\begin{equation}
T(\lambda)=\frac{S_3(\lambda\zeta_2)-\zeta_2^r}{S_2(\lambda)\zeta_2^r-1};\quad U(\lambda)=\frac1{\sqrt3}\zeta_3C(\lambda)\frac{1-S_2^{-1}(\lambda)}{S_3(\lambda\zeta_2)-1}.\label{eq3.32}
\end{equation}
\end{remark}

Elements of the set $E_\alpha$ \eqref{eq2.35} are situated on the bundle of straight lines $L=\bigcup\limits_kL_{\zeta_k}$.

\begin{picture}(200,200)
\put(0,100){\vector(1,0){200}}
\put(150,0){\vector(-1,2){100}}
\put(150,200){\vector(-1,-2){100}}
\put(30,190){$L_{\zeta_2}$}
\put(75,190){$S_2$}
\put(135,150){$\zeta_3w_s$}
\put(175,150){$S_1$}
\put(124,150){$\circ$}
\put(150,96){$\times$}
\put(160,105){$z_k$}
\put(20,150){$S_3$}
\put(20,50){$S_4$}
\put(72,50){$\times$}
\put(45,96){$\circ$}
\put(55,105){$w_s$}
\put(70,147){$\times$}
\put(75,157){$\zeta_2z_k$}
\put(90,20){$S_5$}
\put(105,103){$0$}
\put(190,105){$L_{\zeta_1}$}
\put(123,45){$\circ$}
\put(130,50){$\zeta_2w_s$}
\put(170,60){$S_6$}
\put(80,40){$\zeta_3z_k$}
\put(60,0){$L_{\zeta_3}$}
\end{picture}

\hspace{20mm} Fig. 5

Enumerate them in the following way. Points $\mu_k$ from $E_\alpha$ lying on the axis $L_{\zeta_1}$ enumerate in ascending order, besides, positive $\mu_k$ denote by $\{z_k\}_1^p$ and the negative by $\{w_s\}_1^l$ ($l+p=n<\infty$). Other elements of $E_\alpha$ are obtained from $\{z_k\}_1^p$ and $\{w_s\}_1^l$ after rotations $\lambda\rightarrow\lambda\zeta_2$ and $\lambda\rightarrow\lambda\zeta_3$ (see Fig.5). So, e. g., elements of $E_\alpha$ belonging to the straight line $L_{\zeta_2}$ are given by $\{\zeta_2z_k\}_1^p$ and $\{z_2w_s\}_1^l$. Thus,
\begin{equation}
\begin{array}{ccc}
E_\alpha=\{\zeta_2^lz_k,\zeta_2^lw_s:\{z_k\}_1^p\cup\{w_s\}_1^l\in L_{\zeta_1};l=0,1,2\,(1\leq k\leq p,1\leq s\leq l,\\
p+l=n<\infty)\}
\end{array}\label{eq3.33}
\end{equation}

Equation \eqref{eq3.26} implies that the function
\begin{equation}
y_1(\lambda_m,x)=(0,\psi_1(\lambda_m,x))\label{eq3.33}
\end{equation}
where $\lambda_m\in\left(\{\zeta_2z_k\}_1^p\cup\{\zeta_3w_s\}_1^l\right)$ ($\subset\overline{S_2}$) and
\begin{equation}
\psi_1(\lambda_m,x)=\frac{\alpha i}{\lambda_m^2}\widetilde{q}_1^m(\lambda_m)\int\limits_x^\infty s_2(i\lambda_m(x-t))q(t)dt\label{eq3.34}
\end{equation}
(here $b(\lambda_m)=0$) belongs to $\mathfrak{D}(\mathcal{L}_\alpha)$ \eqref{eq2.4} and is an eigenfunction of the operator $\mathcal{L}_\alpha$, $\mathcal{L}_\alpha y_1(\lambda_m,x)=\lambda_m^3y_1(\lambda_m,x)$. Thus, $\{y_1(\lambda_m,x)\}$ form a set of connected states \cite{6,1,25,3} of the operator $\mathcal{L}_\alpha$ \eqref{eq2.2} where $\lambda_m$ are situated on the boundary of the sector $S_2$.

Similarly, Jost solutions \eqref{eq2.27}
\begin{equation}
y_k(\lambda,x)=(a(\lambda\zeta_k)e^{i\lambda^3x},\psi_k(\lambda,x))\quad(k=2,3)\label{eq3.35}
\end{equation}
obtained from $y_1(\lambda,x)$ \eqref{eq3.26} via substitutions $\lambda\rightarrow\lambda\zeta_2$, $\lambda\rightarrow\lambda\zeta_3$ give two additional series of connected states enumerated by the numbers $(\{z_k\}_1^p\cup\{\zeta_2w_s\}_1^s)\subset\overline{S}_6$ (for $y_2(\lambda,x)$) and numbers $(\{\zeta_3z_k\}_1^p\cup\{w_s\}_1^s)\subset\overline{S}_4$ (for $y_3(\lambda,x)$). These eigenfunctions are easily calculated from
$y_1(\lambda_m,x)$ \eqref{eq3.35} after changes of variables $\lambda_m\rightarrow\lambda_m\zeta_2$, $\lambda_m\rightarrow\lambda_m\zeta_3$.

\section{Inverse problem}

{\bf 4.1.} One-dimensional perturbation $\alpha\langle u,q\rangle q$ in \eqref{eq2.2} is invariant relative to transformation $(\alpha,q)\rightarrow(\alpha c^{-2},qc)$ where $c$ is a real non-zero constant. Therefore the perturbation $\alpha\langle u,q\rangle q$ is in need of normalization. Hereinafter, consider real-valued functions $q(x)$ from $L^2(\mathbb{R}_+)$ for which \eqref{eq2.3} takes place and
\begin{equation}
\left|\int\limits_0^\infty xq(x)dx\right|^2=1.\label{eq4.1}
\end{equation}

Proceed to the step-by-step procedure of recovery of a number $\alpha\in\mathbb{R}$ and a real-valued function $q(x)$ normalized by condition \eqref{eq4.1} from the scattering coefficients $S_2(\lambda)$, $S_3(\lambda)$ and a function $C(\lambda)$ such that $C(\lambda\zeta_2)=\zeta_2^rC(\lambda)$.

{\bf Step 1}. By the scattering coefficients $S_2(\lambda)$, $S_3(\lambda)$ and a function $C(\lambda)$ such that $C(\lambda\zeta_2)=\zeta_2^rC(\lambda)$ holomorphic in $\mathbb{D}_{a'}$ ($a'=a/3$) construct the functions \eqref{eq3.32}
\begin{equation}
T(\lambda)=\frac{S_3(\lambda\zeta_2)-\zeta_2^r}{S_2(\lambda)\zeta_2^r-1};\quad U(\lambda)=\frac1{\sqrt3}\zeta_3C(\lambda)\frac{1-S_2^{-1}(\lambda)}{S_3(\lambda\zeta_2)-1}\label{eq4.2}
\end{equation}
satisfying relations \eqref{eq3.13} and \eqref{eq3.16}.

{\bf Step 2.} Using \eqref{eq3.26}, by $T(\lambda)$ and $U(\lambda)$ \eqref{eq4.2} define the functions $\{G_k(\lambda)\}_1^3$ on the system of rays $l=\bigcup\limits_kl_{\zeta_k}$ \eqref{eq1.16}. Using these functions, on the bundle of straight lines $L=\bigcup\limits_kL_{\zeta_k}$ \eqref{eq1.9} define the functions $\{\widehat{G}_k(\lambda)\}_1^3$,
\begin{equation}
\widehat{G}_k(\lambda)=G_k(\lambda)\chi_{l_{\zeta_k}}+G_k(\lambda\zeta_2)\chi_{\widehat{l}_{\zeta_k}}\quad(1\leq k\leq3)\label{eq4.3}
\end{equation}
where $\chi_{l_{\zeta_k}}$ and $\chi_{\widehat{l}_{\zeta_k}}$ are characteristic functions of the rays $l_{\zeta_k}$ and $\widehat{l}_{\zeta_k}$ \eqref{eq1.10}.

{\bf Step 3.} The functions $\{\widehat{G}_k(\lambda)\}_1^3$ \eqref{eq4.3} are coefficients of the boundary value Riemann problem \eqref{eq3.21} -- \eqref{eq3.23} on the complex contour $L=\bigcup\limits_kL_{\zeta_k}$, which is solved in a standard way \cite{26, 27}. So, in the case if the function $\psi_1(\lambda,0)$ ($\psi_2(\lambda,0)$ and $\psi_3(\lambda,0)$) does not vanish on the closure of the set $\overline{S}_2$ (correspondingly, $\overline{S}_6$ and $\overline{S}_4$), and thus $\psi_1^*(\lambda,0)$, $\psi_2^*(\lambda,0)$, $\psi_3^*(\lambda,0)$ also don't vanish in $\overline{S_5}$, $\overline{S_3}$, $\overline{S_1}$, i. e., index of the boundary value problem \eqref{eq3.21} -- \eqref{eq3.23} vanishes, then the function
\begin{equation}
\psi(z)\stackrel{\rm def}{=}\exp\left\{\sum\limits_k\frac1{2\pi i}\int\limits_{L_{\zeta_k}}\frac{\ln\widehat{G}_k(\lambda)}{\lambda-z}d\lambda\right\}\label{eq4.4}
\end{equation}
gives solution to the boundary value problem \eqref{eq3.21} -- \eqref{eq3.23} and is normalized by the condition $\psi(\infty)=1$. Thus, in this case by the coefficients $\{G_k(\lambda)\}_1^3$ \eqref{eq3.26} set on $l$ \eqref{eq1.16} the functions $\{\psi_k(\lambda,0)\}_1^3$ and $\{\psi_k^*(\lambda,0)\}$ are defined unambiguously. These functions are holomorphic inside corresponding sectors $S_p$ \eqref{eq1.11}. Analyticity of the coefficients $\{G_k(\lambda)\}_1^3$ for $\lambda\in\mathbb{D}_{a'}$ implies holomorphy of the functions $\{\psi_k(\lambda,0)\}_1^3$ and $\{\psi_k^*(\lambda,0)\}_1^3$ in the same domain.

\begin{remark}\label{r4.1}
The function $\psi(z)$ \eqref{eq4.4} equals to the product
\begin{equation}
\psi(z)=\Phi_1(\lambda)\cdot\Phi_2(\lambda)\cdot\Phi_3(\lambda)\label{eq4.5}
\end{equation}
where
\begin{equation}
\Phi_k(\lambda)\stackrel{\rm def}{=}\exp\left\{\frac1{2\pi i}\int\limits_{L_{\zeta_k}}\frac{\ln\widehat{G}_k(\lambda)}{\lambda-z}d\lambda\right\}\quad(1\leq k\leq3)\label{eq4.6}
\end{equation}
and every function $\Phi_k(z)$ is a solution to the Riemann problem on the straight line $L_{\zeta_k}$,
\begin{equation}
\Phi_k(\lambda+i0)=\widehat{G}_k(\lambda)\Phi_k(\lambda-i0)\quad(\lambda\in L_{\zeta_k},1\leq k\leq3),\label{eq4.7}
\end{equation}
besides,
$$\Phi(\lambda\pm i0)=\lim\limits_{\varepsilon\rightarrow+0}\Phi(\zeta_k(\mu\pm i\varepsilon))\quad(\lambda=\mu\zeta_k\in L_{\zeta_k},\mu\in\mathbb{R},\,1\leq k\leq3).$$
\end{remark}

{\bf Step 4.} By the functions $\psi_2^*(\lambda,0)$, $\psi_3^*(\lambda,0)$ and $\psi_1(\lambda,0)$ obtained at the previous step, one unambiguously defines the function, according to formula \eqref{eq3.21},
\begin{equation}
M(\lambda)\stackrel{\rm def}{=}-\frac{\alpha i}{3\lambda^2}(\widetilde{q}_2(\lambda)-\widetilde{q}_3(\lambda))(\widetilde{q}_2(\lambda)-\widetilde{q}_3(\lambda))^*\quad(\forall\lambda\in\mathbb{D}_{a'}).\label{eq4.8}
\end{equation}
Since
$$M(0)=-\alpha i\left|\int\limits_0^\infty xq(x)dx\right|^2,$$
then, using normalization \eqref{eq4.1}, one finds $\alpha=iM(0)$. Therefore from \eqref{eq4.8} one uniquely defines the function
\begin{equation}
N(\lambda)=Q(\lambda)Q^*(\lambda);\quad Q(\lambda)\stackrel{\rm def}{=}\widetilde{q}_2(\lambda)-\widetilde{q}_3(\lambda).\label{eq4.9}
\end{equation}
Reality of $q(x)$ implies that $Q^*(\lambda)=-Q(-\overline{\lambda})$ and thus
\begin{equation}
N(\lambda)=-Q(\lambda)Q(-\lambda)\quad(\forall\lambda\in\mathbb{D}_{a'}\cap\mathbb{R}),\label{eq4.10}
\end{equation}
therefore
\begin{equation}
N(-\lambda)=N(\lambda),\quad N(\lambda)=|Q(\lambda)|^2\quad(\forall\lambda\in\mathbb{D}_{a'}\cap\mathbb{R}).\label{eq4.11}
\end{equation}
By the holomorphic in $\mathbb{D}_{a/3}$ function $N(\lambda)$ from \eqref{eq4.10} unambiguously calculate an analytic for all $\lambda\in\mathbb{D}_{a'}\cap\mathbb{R}$ function $|Q(\lambda)|=\sqrt{N(\lambda)}$ for which the equality $|Q(-\lambda)|=|Q(\lambda)|$ holds for all $\lambda\in\mathbb{D}_{a'}\cap\mathbb{R}$ (see \eqref{eq4.11}). Consider a real-valued analytic function $\varphi(\lambda)$ defined for $\lambda\in\mathbb{D}_{a/3}$ such that
$$\varphi(\lambda)+\varphi(-\lambda)=\pi\quad(\forall\lambda\in\mathbb{D}_{a'}\cap\mathbb{R}).$$
Define the function
\begin{equation}
Q(\lambda)\stackrel{\rm def}{=}\sqrt{N(\lambda)}e^{i\varphi(\lambda)}\quad(\lambda\in\mathbb{D}_{a'}\cap\mathbb{R})\label{eq4.12}
\end{equation}
analytic for all $\lambda\in\mathbb{D}_{a'}\cap\mathbb{R}$ such that equalities $Q(-\lambda)=-Q(\lambda)$ and $N(\lambda)=|Q(\lambda)|^2$ hold for the same values of $\lambda$.

So, at this step, one calculates the number $\alpha\in\mathbb{R}$ and recovers (ambiguously) the function $Q(\lambda)$ \eqref{eq4.9} from $N(\lambda)$ \eqref{eq4.9} that is naturally holomorphically extendable into $\mathbb{D}_{a'}$.

{\bf Step 5.} Show that by $Q(\lambda)$ \eqref{eq4.9} the function $q(x)$ is found unambiguously. Obviously, the function
\begin{equation}
Q(\lambda)=\widetilde{q}_2(\lambda)-\widetilde{q}_3(\lambda)=\int\limits_0^\infty\left(e^{-i\lambda\zeta_2x}-e^{-i\lambda\zeta_3x}\right)q(x)dx\label{eq4.13}
\end{equation}
is holomorphically extendable into the sector $S_2$ and its values on the imaginary axis $i\mathbb{R}_+$ ($\lambda=i\mu$, $\mu\geq0$) equal
\begin{equation}
Q(i\mu)=2i\int\limits_0^\infty e^{-\frac\mu 2}\sin\frac{\sqrt3\mu x}2q(x)dx.\label{eq4.14}
\end{equation}
So, one knows the function
\begin{equation}
F(y)=\int\limits_0^\infty e^{-yt}\sin\sqrt3 ytq(t)dt\quad(\mu=2y\geq0)\label{eq4.15}
\end{equation}
which is a hybrid of two well-known integral transforms, Laplace transform and Fourier sine-transform. Inverse of the transformation $q(x)\rightarrow F(y)$ \eqref{eq4.15} as well as description of image and preimage of this correspondence is unknown.

Give method of recovery of the function $q(x)$ from $F(y)$ \eqref{eq4.15}, to do this consider the function of two real variables
\begin{equation}
\Phi(x,y)\stackrel{\rm def}{=}\int\limits_0^\infty e^{-yt}\sin xtq(t)dt\label{eq4.16}
\end{equation}
given in the half-plane $(x,y)\in\mathbb{R}\times\mathbb{R}_+$. The function $\Phi(x,y)$ \eqref{eq4.16} is harmonic in this half-plane, $\Delta\Phi=0$, and its boundary value is
\begin{equation}
\Phi(x,0)=\int\limits_0^\infty\sin xtq(t)dt\stackrel{\rm def}{=}g(x),\label{eq4.17}
\end{equation}
besides, $g(-x)=-g(x)$. The function $\Phi(x,y)$ \eqref{eq4.15} is expressed via its boundary values $g(x)$ \eqref{eq4.17} using the Poisson kernel \cite{28, 24},
\begin{equation}
\Phi(x,y)=\frac1\pi\int\limits_{\mathbb{R}}\frac y{(x-t)^2+y^2}g(t)dt,\label{eq4.18}
\end{equation}
besides, $\Phi(-x,y)=-\Phi(x,y)$. The function $F(y)=\Phi(y\sqrt{3},y)$ is known.

\begin{picture}(300,100)
\put(0,10){\vector(1,0){300}}
\put(150,0){\vector(0,1){100}}
\put(150,10){\line(3,1){90}}
\put(150,10){\line(-3,1){90}}
\put(155,0){$0$}
\put(290,12){$x$}
\put(153,95){$y$}
\put(75,40){$x=-\sqrt3 y$}
\put(0,50){$-\Phi(y\sqrt3,y)$}
\put(155,40){$\sqrt3 y=x$}
\put(220,50){$\Phi(y\sqrt3,y)$}
\end{picture}

\hspace{20mm} Fig. 6

Extend $F(y)$ onto a straight line $x=-\sqrt3y$ by symmetry $\widetilde{F}(y)=-\Phi(y\sqrt3,y)$. The function $\Phi(x,y)$ is harmonic inside the angle from upper half-plane, sides of which are given by the straight lines $x=\pm\sqrt3y$ (Fig. 6) and is known on its boundary (on the sides of the angle). Consequently, as a harmonic function, it is defined by its boundary values $F(y)$ and $\widetilde{F}(y)$. Set $x=0$ in equation \eqref{eq4.18}, then
$$\Phi(0,y)=\frac1\pi\int\limits_0^\infty\frac y{t^2+y^2}g(t)dt\quad(y\in\mathbb{R}_+).$$
Use
$$\int\limits_0^\infty\frac y{y^2+a^2}\sin xydy=\frac\pi2e^{ax}\quad(a>0),$$
hence one finds that
\begin{equation}
G(y)\stackrel{\rm def}{=}\int\limits_0^\infty\Phi(0,y)\sin yxdy=\frac12\int\limits_0^\infty e^{-tx}g(t)dt\quad(x>0).\label{eq4.19}
\end{equation}
So, the Laplace transform $G(y)$ \eqref{eq4.19} of the function $g(x)$ \eqref{eq4.17} is known. The function $G(y)$ is holomorphically extendable up to $G(z)$ where $z$ lies in the right half-plane, $\Re z>0$. Using the inverse Laplace transform, calculate $g(x)$,
$$g(x)=\frac1{2\pi i}\int\limits_{c-i\infty}^{c+i\infty}e^{zt}G(z)dz.$$
Since $g(x)$ is given by \eqref{eq4.17}, then using the inverse Fourier sine-transform one recovers the function $q(x)$.

So, the function $Q(\lambda)$ can be represented as a difference of Fourier transforms \eqref{eq4.13} and is holomorphic in the sector $S_2$. Analyticity of $Q(\lambda)$ in the region $\mathbb{D}_{a'}$ implies analyticity of the functions $\widetilde{q}_2(\lambda)$ and $\widetilde{q}_2(\lambda)$ in the same domain, which due to Paley -- Wiener theorem \cite{23,24} leads to condition \eqref{eq2.3} for a function $q(x)$ for some $a$ ($>0$).

{\bf Conclusion 1.} {\it By holomorphic in the region $\mathbb{D}_{a'}$ ($a'>0$) functions $S_2(\lambda)$, $S_3(\lambda)$ and $C(\lambda)$ where $C(\lambda\zeta_2)=\zeta_2^rC(\lambda)$ ($r=0$, $1$, $2$) satisfying the conditions \eqref{eq3.13}, \eqref{eq3.16} the number $\alpha\in\mathbb{R}$ and the real-valued function $q(x)$ normalized by condition \eqref{eq4.1}, for which \eqref{eq2.3} takes place for some $a>0$, are recovered.}
\vspace{5mm}

{\bf 4.2.} Describe the method of recovery of the number $\alpha$ and the function $q(x)$ in the presence of connected states to which there corresponds the set $E_\alpha$ \eqref{eq3.33} situated on the bundle of straight lines $L=\bigcup\limits_k L_{\zeta_k}$ (see Fig. 7).

\begin{picture}(200,200)
\put(0,100){\vector(1,0){200}}
\put(150,0){\vector(-1,2){100}}
\put(150,200){\vector(-1,-2){100}}
\put(30,190){$L_{\zeta_2}$}
\put(95,190){$S_2$}
\put(100,150){$\centerdot$}
\put (100,159){$\theta_2$}
\put(175,150){$S_1$}
\put(150,130){$\centerdot$}
\put(160,130){$\theta_1$}
\put(150,70){$\centerdot$}
\put(160,60){$\theta_6$}
\put(100,40){$\centerdot$}
\put(20,150){$S_3$}
\put(20,50){$S_4$}
\put(90,20){$S_5$}
\put(100,50){$\theta_5$}
\put(40,60){$\centerdot$}
\put(50,60){$\theta_4$}
\put(40,140){$\centerdot$}
\put(50,140){$\theta_3$}
\put(190,105){$L_{\zeta_1}$}
\put(170,40){$S_6$}
\put(60,0){$L_{\zeta_3}$}
\end{picture}

\hspace{20mm} Fig. 7

Choose on the bisector of the sector $S_1$ a point $\theta_1$ and let
\begin{equation}
\theta_k=\eta^{k-1}\theta_1\quad(1\leq k\leq6)\label{eq4.20}
\end{equation}
where $\eta=e^{i\pi/6}$. Points $\{\theta_k\}_1^6$ are situated on bisectors of corresponding sectors $\{S_k\}_1^6$ \eqref{eq1.11} (see Fig. 7). Write the boundary value problem \eqref{eq3.21} as
\begin{equation}
\begin{array}{ccc}
{\displaystyle\left(\frac{\lambda-\theta_2}{\lambda-\theta_3}\right)^p\cdot G_2(\lambda)\cdot\psi_1(\lambda,0)\cdot\frac{(\lambda-\theta_3)^p}{\displaystyle\prod\limits_{k=1}^p(\lambda-\zeta_2z_k)}=\psi_2^*(\lambda,0)\cdot\frac{(\lambda-\theta_2)^p}{\displaystyle\prod\limits_{k=1}^p(\lambda-\zeta_2z_k)}
\quad(\lambda\in l_{\zeta_2});}\\
{\displaystyle\left(\frac{\lambda-\theta_6}{\lambda-\theta_5}\right)^l\cdot G_2(\lambda\zeta_2)\cdot\psi_2(\lambda,0)\cdot\frac{(\lambda-\theta_5)^l}{\displaystyle\prod\limits_{s=1}^l(\lambda-\zeta_2w_s)}=\psi_1^*(\lambda,0)\frac{(\lambda-\theta_6)^l}{\displaystyle\prod\limits_{s=1}^l(\lambda-\zeta_2w_s)}
\quad(\lambda\in\widehat{l}_{\zeta_2});}
\end{array}\label{eq4.21}
\end{equation}
problem \eqref{eq3.22} write as
\begin{equation}
\begin{array}{ccc}
{\displaystyle\left(\frac{\lambda-\theta_4}{\lambda-\theta_5}\right)^p\cdot G_3(\lambda)\cdot\psi_3(\lambda,0)\frac{(\lambda-\theta_5)^p}{\displaystyle\prod\limits_{k=1}^p(\lambda-\zeta_3z_k)}=\psi_1^*(\lambda,0)\frac{(\lambda-\theta_4)^p}{\displaystyle\prod\limits_{k=1}^p(\lambda-\zeta_3z_k)}\quad(l\in l_{\zeta_3})}\\
{\displaystyle\left(\frac{\lambda-\theta_2}{\lambda-\theta_1}\right)^l\cdot G_3(\lambda\zeta_2)\cdot\psi_1(\lambda,0)\frac{(\lambda-\theta_1)^l}{\displaystyle\prod\limits_{s=1}^l(\lambda-\zeta_3w_s)}=\psi_3^*(\lambda,0)\frac{(\lambda-\theta_2)^l}{\displaystyle\prod\limits_{s=1}^l(\lambda-\zeta_3w_s)}\quad
(l\in\widehat{l}_{\zeta_3});}
\end{array}\label{eq4.22}
\end{equation}
and, finally, write problem \eqref{eq3.23} as
\begin{equation}
\begin{array}{ccc}
{\displaystyle\left(\frac{\lambda-\theta_6}{\lambda-\theta_1}\right)^pG_1(\lambda)\psi_2(\lambda,0)\frac{(\lambda-\theta_1)^p}{\displaystyle\prod\limits_{k=1}^p(\lambda-z_k)}=\psi_3^*(\lambda,0)\frac{(\lambda-\theta_6)^p}
{\displaystyle\prod\limits_{k=1}^p(\lambda-z_k)}\quad(\lambda\in l_{\zeta_1})};\\
{\displaystyle\left(\frac{\lambda-\theta_4}{\lambda-\theta_3}\right)^lG_1(\lambda\zeta_2)\psi_3(\lambda,0)\frac{(\lambda-\theta_3)^s}{\displaystyle\prod\limits_{s=1}^l(\lambda-w_s)}=\psi_2^*(\lambda,0)\frac{(\lambda-\theta_4)^l}
{\displaystyle\prod\limits_{s=1}^l(\lambda-w_s)}\quad(\lambda\in\widehat l_{\zeta_1}).}
\end{array}\label{eq4.23}
\end{equation}
Notice that the function ${\displaystyle\psi_1(\lambda,0)(\lambda-\theta_3)^p\left(\prod\limits_{k=1}^p(\lambda-\zeta_2z_k)\right)^{-1}}$ from the boundary value problem \eqref{eq4.21} does not vanish on $l_{\zeta_2}$ and tends to identity as $\lambda\rightarrow\infty$. Other similar functions from the boundary value problems \eqref{eq4.21} -- \eqref{eq4.23} have the same properties. Therefore, using technique of Subsec. 4.1, one recovers the number $\alpha\in\mathbb{R}$ and real-valued function $q(x)$. As always (see \cite{26,27}), solution to the boundary value problem \eqref{eq4.21} -- \eqref{eq4.23} does not depend on the choice of $\theta$.
\vspace{5mm}

{\bf 4.3.} The totality
\begin{equation}
\Omega=\{S_2(\lambda);S_3(\lambda);C(\lambda);E_\alpha\}\label{eq4.24}
\end{equation}
is a scattering problem data. Description of scattering data is given in the theorem.

\begin{theorem}\label{t4.1}
Let scattering data $\Omega$ \eqref{eq4.24} satisfy the requirements:

$\rm{(i)}$ $S_2(\lambda)$, $S_3(\lambda)$, $C(\lambda)$ are given for all $\lambda\in\mathbb{C}$ and are holomorphic for $\lambda\in\mathbb{D}_{a'}=\{\lambda\in\mathbb{C}:|\lambda|<a',a'>0\}$;

$\rm{(ii)}$ the function $c(\lambda)$ has the property $C(\lambda\zeta_2)=\zeta_2^rC(\lambda)$ ($r$ is a fixed number from $\mathbb{N}$) and constructed by $S_2(\lambda)$, $S_3(\lambda)$, $C(\lambda)$ functions $T(\lambda)$ and $U(\lambda)$ \eqref{eq4.2} satisfy the equalities \eqref{eq3.13}, \eqref{eq3.16};

$\rm{(iii)}$ finite set $E_\alpha$ is situated on the bundle of straight lines $L=\bigcup\limits_kL_{\zeta_k}$ \eqref{eq1.23} and is given by \eqref{eq3.33}.

Then there exists an operator $\mathcal{L}_\alpha$ \eqref{eq2.2}, \eqref{eq2.4}, spectral data of which is the set $\Omega$, besides, $\alpha\in\mathbb{R}$, $q(x)$ is a real-valued function from $L^2(\mathbb{R}_+)$ for which for some $a$ \eqref{eq2.3} takes place and normalization \eqref{eq4.1} is true.
\end{theorem}

\renewcommand{\refname}{ \rm 
\end{Large}
\end{document}